# Do Gödel's incompleteness theorems set absolute limits on the ability of the brain to express and communicate mental concepts verifiably?

## A lay perspective[1]

Bhupinder Singh Anand[2]

Classical interpretations of Gödel's formal reasoning, and of his conclusions, implicitly imply that mathematical languages are essentially incomplete, in the sense that the truth of some arithmetical propositions of any formal mathematical language, under any interpretation, is, both, non-algorithmic, and essentially unverifiable. However, a language of general, scientific, discourse, which intends to mathematically express, and unambiguously communicate, intuitive concepts that correspond to scientific investigations, cannot allow its mathematical propositions to be interpreted ambiguously. Such a language must, therefore, define mathematical truth verifiably. We consider a constructive interpretation of classical, Tarskian, truth, and of Gödel's reasoning, under which any formal system of Peano Arithmetic - classically accepted as the foundation of all our mathematical languages - is verifiably complete in the above sense. We show how some paradoxical concepts of Quantum mechanics can, then, be expressed, and interpreted, naturally under a constructive definition of mathematical truth.

---

[2] The author is an independent scholar. E-mail: anandb@vsnl.com. Address: 32, Agarwal House, D Road, Churchgate, Mumbai - 400 020, INDIA. Tel: +91 (22) 2281 3353. Fax: +91 (22) 2209 5091.

# Contents[3]

---

# I Do Gödel's Theorems limit our ability to mathematically express, and communicate, mental concepts?

## I-1 Introduction[4]

Classical[5] interpretations of Gödel's formal reasoning and conclusions - in his seminal 1931 paper, "*On formally undecidable propositions of Principia Mathematica and related systems I*" [Go31a], in which he introduces his two, famous, "Incompleteness" Theorems - suggest that the classical, Tarskian, truth[6] of some propositions[7] of a formal mathematical language[8], under an interpretation[9], is, both, non-algorithmic[10], and essentially unverifiable, constructively[11].

---

[4] We define the notation and terminology used in this essay in the Appendix.

[5] For the purposes of this essay, we take the expositions by Hardy [Ha47], Landau [La51], Mendelson [Me64], Rudin [Ru53] and Titchmarsh [Ti61] as standard presentations of classical mathematical reasoning and conclusions.

[6] We use the word "true" both in its familiar, linguistic, sense, and in a mathematically precise sense; the appropriate meaning is usually obvious from the context. Mathematically, we follow Mendelson's exposition of the truth of a formal sentence under an interpretation as determined by Tarski's definitions of satisfiability and truth ([Me64], p51).

[7] When referring to a formal language, we assume the terms "sentence" and "proposition" are synonymous, and that they refer to a well-formed expression of the language that contains no free variables, and which translates, under an interpretation, as a proposition in the usual, linguistic, sense. However, the term "proposition" is sometimes reserved for "metatheorems", which are proven statements about the language ([Me64], footnote on p31-31).

[8] By a "formal language" we mean a "formal system" or a "formal theory" as in Mendelson ([Me64], p29).

[9] As in this sentence, the word "interpretation" may be used both in its familiar, linguistic, sense, and in a mathematically precise sense; the appropriate meaning is usually obvious from the context. Mathematically, we follow Mendelson's definition of "interpretation" ([Me64], §2, p49):

"An *interpretation* consists of a non-empty set D, called the *domain* of the interpretation, and an assignment to each predicate letter $A_j^n$ of an $n$-place relation in D, to each function letter $f_j^n$ of an $n$-place operation in D (i.e., a function from $D^n$ into D), and to each individual constant $a_i$ of some fixed element of D. Given such an interpretation, variables are thought of as ranging over the set D, and ~, =>, and quantifiers are given their usual meaning. (Remember that an $n$-place relation in D can be thought of as a subset of $D^n$, the set of all $n$-tuples of elements of D.)"

The questions arise: Does this imply, first, that the determination of mathematical truths, akin to that of scientific truths, is a process best described as discovery; and, second, do Gödel's Theorems set absolute limits on our ability to mathematically express, and communicate, our mental concepts precisely, and verifiably, in general, scientific, discourse[12]?

Since it is not germane to the issue in these essays, we shall neither divert ourselves with an exposition of the novel meta-mathematical, and logical, features of these remarkable

---

We note that the interpreted relation $R(x)$ is obtained from the formula $[R(x)]$ of a formal system P by replacing every primitive, undefined symbol of P in the formula $[R(x)]$ by an interpreted mathematical symbol (i.e. a symbol that is a shorthand notation for some, semantically well-defined, concept of classical mathematics). So the P-formula $[(Ax)R(x)]$ interprets as the sentence $(Ax)R(x)$, and the P-formula $[\sim(Ax)R(x)]$ as the sentence $\sim(Ax)R(x)$.

We also note that the meta-assertions "$[(Ax)R(x)]$ is a true sentence under the interpretation M of P", and "$(Ax)R(x)$ is a true sentence of the interpretation M of P", are equivalent to the meta-assertion "$R(x)$ is satisfied for any given value of $x$ in the domain of the interpretation M of P" ([Me64], p51).

[10] We follow Mendelson's definition of an algorithm as an effectively computable function ([Me64], p208), with the proviso that such computability is, necessarily, by some uniform rule, or method. Subject to their being individually proven as formal mathematical objects - a concept that we define precisely - we also follow Mendelson's set-theoretic definitions of a "function" and of a "relation" ([Me64], p6-7). In other words, since we intend to argue later that there may be number-theoretic functions that are not definable as formal mathematical objects in any Axiomatic Set Theory, we treat the sets in Mendelson's definitions as hypothetical and intuitive, but not formal, mathematical objects. Therefore, assuming formal set-theoretic properties for them, even in informal reasoning, may invite inconsistency.

[11] The term "constructive" is used both in its familiar, linguistic, sense, and in a mathematically precise sense. Mathematically, we term a concept as "constructive" if, and only if, it can be defined in terms of pre-existing concepts without inviting inconsistency (cf. Mendelson's remarks in [Me64], p82). Otherwise, we understand it in an intuitive sense to mean unambiguously verifiable, by some "effective method" ([Me64], p207-8), within some finite, well-defined, language or meta-language ([Me64], p31, footnote). Generally, it may be taken to correspond, broadly, to Gödel's concept of "intuitionistically unobjectionable" ([Go31a], p26).

However, the concept of "effective method" (as well its synonym, "mechanical procedure") are not at all precise in classical theory. As Mendelson notes ([Me64], p207), "... what we mean is a process which requires no ingenuity for its performance". In these essays, we consider some consequences of defining these critical concepts precisely in appropriate contexts.

[12] We briefly discuss some of the issues involved in the representation of mental concepts within a formal language in §II-5, "Can all mental concepts be expressed mathematically?".

meta-theorems, nor with details of their fascinating meta-proofs[13]. For the purist, Gödel's original, 1931, paper (cf. [Go31a] and [Go31b]) remains unsurpassed as the definitive source for a study of the Theorems, and of their immediate consequences. For the interested, popular discussions of the Theorems, and of their commonly perceived meanings and implications, are lucidly provided, and extensively referenced, by Penrose in [Pe90], and in [Pe94].

What interests us here, rather, is the possibility that there may be a "loophole", in the classical interpretation of Gödel's conclusion, which allows us to define mathematical truth in an effectively verifiable way. We, thus, limit ourselves to reviewing, essentially from a layperson's perspective, the argument that any formal system of Peano Arithmetic[14], PA, based on Dedekind's formulation of the Peano Postulates[15], is essentially incomplete[16].

---

[13] We use words such as "proof", "meta-mathematics", "meta-theorem", "meta-proof", "meta-language", etc., both in their familiar linguistic sense, and in a mathematically precise sense; the appropriate meaning is usually obvious from the context. As Mendelson notes (cf. [Me64], p29 and footnote on p31):

"The word 'proof' is used in two distinct senses. First, it has a precise meaning defined above as a certain finite sequence of wfs of L. However, in another sense, it also designates certain sequences of sentences of the English language (supplemented by various technical terms) which are supposed to serve as an argument justifying some some assertions about the language L (or other formal theories). In general, the language we are studying (in this case L) is called the *object language*, while the language in which we formulate and prove results about the object language is called the *metalanguage*. The metalanguage might also be formalized and made the subject of study, which we would carry out in a meta-metalanguage, etc. However, we shall use the English language as our (unformalized) metalanguage, although, for a substantial part of this book, we employ only a mathematically weak portion of the English language. The contrast between object language and metalanguage is also present in the study of a foreign language; for example, in a German class, German is the object language, while the metalanguage, the language we use, is English. The distinction between 'proof' and 'metaproof' (i.e., a proof in the metalanguage) leads to a distinction between theorems of the object language and *metatheorems* of the metalanguage. To avoid confusion, we generally use 'proposition' instead of 'metatheorem'. The word 'metamathematics' refers to the study of logical and mathematical object languages; sometimes the word is restricted to those investigations which use what appear to the metamathematician to be constructive (or so-called *finitary*) methods."

[14] By a formal system of Peano Arithmetic, we mean an intended formalisation, of Dedekind's formulation of the intuitively-interpreted Peano Axioms, such as Gödel's formal system P ([Go31a], p9), Mendelson's first order theory S ([Me64], p102), or Podnieks first order arithmetic PA2 ([Po01], §3.1).

We note that Mendelson's first-order theory S:

Specifically, we consider whether the classical interpretation of Gödel's argument - that there is some undecidable[17] arithmetical proposition (which we may refer to as GUS, and formally write as the well-formed formula[18] $[(\mathrm{A}x)R(x)]$) that is unprovable[19] in PA, but

---

"...has a single predicate letter $A_1^2$ (as usual, we write $t = s$ for $A_1^2$ $(t, s)$; it has one individual constant $a_1$ (written, as usual, 0); it has three function letters $f_1^1, f_1^2, f_2^2$. We shall write $t'$ instead of $f_1^1(t)$; $t + s$ instead of $f_1^2(t, s)$; and $t*s$ instead of $f_2^2(t, s)$", and the following axioms:

(S1) $(x_1 = x_2) \Rightarrow ((x_1 = x_3) \Rightarrow (x_2 = x_3))$
(S2) $(x_1 = x_2) \Rightarrow (x_1' = x_2')$
(S3) $0 =/= x_1'$
(S4) $(x_1' = x_2') \Rightarrow (x_1 = x_2)$
(S5) $(x_1 + 0) = x_1$
(S6) $(x_1 + x_2') = (x_1 + x_2)'$
(S7) $(x_1 * 0) = 0$
(S8) $(x_1 * x_2') = ((x_1 * x_2) + x_1)$
(S9) For any well-formed formula $F(x)$ of S, $F(0) \Rightarrow ((\mathrm{A}\ x)(F(x) \Rightarrow F(x')) \Rightarrow (\mathrm{A}\ x)F(x)$"

[15] We follow Mendelson's following formulation ([Me64], p102) of Dedekind's (1901) presentation of the Peano's Postulates:

(P1) 0 is a natural number.
(P2) If $x$ is a natural number, there is another natural number denoted by $x'$ (and called the successor of $x$).
(P3) $0 =/= x'$ for any natural number $x$.
(P4) If $x' = y'$, then $x = y$.
(P5) If $Q$ is a property which may, or may not, hold of natural numbers, and if (I) 0 has the property $Q$, and (II) whenever a natural number $x$ has the property Q, then $x'$ has the property $Q$, then all natural numbers have the property $Q$ (Principle of Induction).

[16] We define a language L as "essentially incomplete" if, and only if, L, and every consistent axiomatic extension L' of L, has some undecidable sentence, i.e., some proposition $F$ such that both $[F]$ and $[\sim F]$ are not provable in the language (each extension may have a different undecidable proposition). We define a language as axiomatic if, and only if, there is an effective method for determining whether any given well-formed formula of the language is an axiom. We define an extension L' of a language L as the language obtained by adding a finite number of well-formed formulas to the axioms of L.

[17] A formal proposition $[R]$ is undecidable in a language L if, and only if, both $[R]$ and $[\sim R]$ are unprovable in L, where, under the standard interpretation of L, $\sim R$ interprets as the negation of $R$.

[18] By a "well-formed formula", we refer to any finite concatenation of the symbols of a formal language, constructed according to specified rules of the grammar of the language for the formation of well-formed formulas, to which we attach no meaning. Where the intention is clear from the context, we may refer to a "well-formed formula" simply as a "formula".

[19] We define a formula $[F]$ of a finitely axiomatisable formal language L as provable in L if, and only if, there is a finite sequence of L-formulas, ending in $[F]$, such that each formula of the sequence is either an

true, under the standard interpretation[20] M of PA - is constructive, and intuitionistically unobjectionable[21].

We then consider whether, under a suitably constructive interpretation of classical, Tarskian, truth, and of Gödel's reasoning, any formal system of Peano Arithmetic - classically accepted as the foundation of all our mathematical languages - is verifiably complete, and indicate some consequences of such an interpretation.

We, finally, address, and review, the questions[22]:

    (1) Are Platonism and Formalism incompatible doctrines?

    (2) Is mathematical truth verifiable effectively?

---

axiom of L, or an immediate consequence of some of the formulas preceding it, in the sequence, by means of the rules of inference of L. We define [$F$] as unprovable in L if, and only if, there is no such sequence.

[20] We follow Mendelson's definitions of "standard interpretation" ([Me64], p107) as:

"... the interpretation in which

($a$) the set of non-negative integers is the domain,
($b$) the integer 0 is the interpretation of the symbol 0,
($c$) the successor operation (addition of 1) is the interpretation of the ' function (i.e., of $f_1^1$),
($d$) ordinary addition and multiplication are the interpretations of + and .,
($e$) the interpretation of the predicate letter = is the identity relation."

In other words, the interpreted relation $R(x)$ is obtained from the formula [$R(x)$] of a formal system P by replacing every primitive, undefined, symbol of P in the formula [$R(x)$] by an interpreted mathematical symbol (i.e. a symbol that is a shorthand notation for some, semantically well-defined, concept of classical mathematics) as in ($a$)-($e$).

[21] We generally use the words "constructive / non-constructive, and intuitionistically unobjectionable / objectionable" in Gödel's, rather broad, sense, as expressed by him at the end of his proof of Theorem VI in his 1931 paper [Go31a] on formally undecidable propositions.

Gödel remarks: "One can easily convince oneself that the proof we have just given is constructive (for all the existential assertions occurring in the proof rest upon Theorem V which, as it is easy to see, is intuitionistically unobjectionable), ...".

[22] Since these issues are sought to be addressed independently, there is an element of repetition - for which the author begs the readers' indulgence - that could, perhaps, have been avoided, but possibly at the expense of readability and clarity of exposition.

(3) What is the significance of Gödel's First Incompleteness Theorem?

(4) What is the significance of Turing's Halting Theorem?

(5) Can all mental concepts be expressed mathematically?

(6) Can Peano Arithmetic model some of the more paradoxical concepts of Quantum Mechanics?

## I-2     Is Gödel's undecidable proposition true in a constructive, and intuitionistically unobjectionable, way?

Now, classically, Gödel's well-formed formula $[(Ax)R(x)]$ does translate as a true sentence under the standard interpretation of PA. However, in general, there is nothing intuitive or constructive - in the sense of being effectively verifiable - about such "truth".

The proposition is "true" only if we accept:

**Tarski's definition**[23]: A well-formed formula $[(Ax)F(x)]$ of a language L is true under an interpretation M of L if, and only if, the interpreted relation[24], $F(x)$, is satisfied by every $x$ in M.

---

[23] We take Mendelson ([Me64], p49-52) as a standard exposition of Tarski's definitions of the "satisfiability" and "truth" of well-formed formulas under a given interpretation:

"The notions of satisfiability and truth are intuitively clear, but, for the skeptical, they can be made precise in the following way (cf. [Ta36], p261-405). Let there be given an interpretation with domain D. Let S be the set of denumerable sequences of elements of D. We shall define what it means for a sequence $s = (b_1, b_2, ...)$ in S to satisfy a wf $A$ under the given interpretation. As a preliminary step we define a function $s*$ of one argument, with terms as arguments and values in D.

(1) If $t$ is $x_i$, let $s*(t)$ be $b_i$.
(2) If $t$ is an individual constant, then $s*(t)$ is the interpretation in D of this constant.
(3) If $f_j^n$ is a function letter and $g$ is the corresponding operation in D, and $t_1, ..., t_n$ are terms, then $s*(f_j^n(t_1, ..., t_n)) = g(s*(t_1), ..., s*(t_n))$.

Although this appears to be a fairly innocent formalisation of intuitive truth, we note, first, that it is silent on the question of how, for any interpretation M, we can effectively determine whether the relation $F(x)$ is, indeed, satisfied by every $x$ in M.

Second, it does not distinguish between languages of expression, which are intended to capture elements, of the mental gestalt of an individual, within a symbolic language (reasonably, these would include our spoken and written languages, as also sign languages, painting, sculpture, music, etc.), and languages of communication, which are

---

Thus, $s*$ is a function, determined by the sequence s, from the set of terms into D. Intuitively, for a sequence $s = (b_1, b_2, ...)$ and a term $t$, $s*(t)$ is the element of D obtained by substituting, for each $i$, $b_i$ for all occurences of $x_i$ in $t$, and then performing the operations of the interpretation corresponding to the function letters of $t$. For instance, if $t$ is $f_2^2(x_3, f_1^2(x_1, a_1))$, and the interpretation has the set of integers as its domain, $f_2^2$ and $f_1^2$ are interpreted as ordinary multiplication and addition, and $a_1$ is interpreted as 2, then, for any sequence $s = (b_1, b_2, ...)$ of integers, $s*(t)$ is the integer $b_3$ x $(b_1 + 2)$.

Now we proceed to the definition proper, which is an inductive definition.

(*i*) If $A$ is an atomic wf $A_j^n(t_1, ..., t_n)$ and $B_j^n$ is the corresponding relation of the interpretation, then the sequence s satisfies $A$ if and only if $B_j^n(s*(t_1), ..., s*(t_n))$, i.e., if the $n$-tuple $(s*(t_1), ..., s*(t_n))$ is in the relation $B_j^n$.

(*ii*) s satisfies ~$A$ if and only if $s$ does not satisfy $A$.

(*iii*) s satisfies $A => B$ if and only if either s does not satisfy $A$ or s satisfies $B$.

(*iv*) s satisfies $(x_i)A$ if and only if every sequence of S which differs from $s$ in at most the $i$'th component satisfies $A$.

Intuitively, a sequence $s = (b_1, b_2, ...)$ satisfies a wf $A$ if and only if, when we substitute, for each $i$, a symbol representing $b_i$ for all free occurences of $x_i$ in $A$, the resulting proposition is true under the given interpretation.

A wf $A$ is *true* (for the given interpretation) if and only if every sequence in S satisfies $A$.

$A$ is *false* (for the given interpretation) if and only if no sequence in S satisfies $A$.

An interpretation is said to be a *model* for a set T of wfs if and only if every wf in T is true for the interpretation."

For a detailed overview of Tarski's truth definitions, see [Hw01].

[24] The word "predicate" is often used as a synonym for "relation".

intended to distinguish, and effectively communicate, those of such individual concepts that are accepted as lying within what may be accepted, and termed, as a common collective of gestalts.

This lack of distinction is reflected in the oft-encountered - and, as we argue, unnecessary[25] - controversy between those who believe that whatever can be conceived must exist in a Platonic world, and those who believe that only that which can be communicated effectively can be claimed to exist.

Moreover, a significant consequence - of the failure to distinguish between Platonic conception and effective communication - is that Tarski's definition implicitly commits us to admitting, in a formal language, such as, say, PA, implicit reference to Platonic elements, in the domain of an interpretation M of PA, that are, clearly, non-intuitive, and conceivable only subjectively in individual gestalts[26].

Thus, the definition (implicitly) implies that we may (explicitly) assert the closure of a formal relation under the universal quantifier as satisfied in M if, and only if, the relation is individually, *and* collectively, satisfied by all the elements in the ontology of M, even if some elements of this ontology are *not* interpretations of any mathematical objects that are representable[27] in PA!

---

[25] We address this issue separately in §II-1, "Are Platonism and Formalism incompatible doctrines?".

[26] In other words, Tarski's definitions implicitly lend legitimacy to the belief that abstract mathematical concepts are objective mathematical realities that can be perceived in a manner analogous to sense perception. Known as mathematical realism, such a philosophy holds that mathematical entities exist independently of the human mind.

The major problem of mathematical realism is this: precisely where and how do the mathematical entities exist? Is there a world, completely separate from our physical one, which is occupied by the mathematical entities? How can we gain access to this separate world and discover truths about the entities?

[27] The words "expressible" and "representable" can, generally, be understood in their familiar, linguistic, sense. However, when used in a mathematical sense (the appropriate sense is usually obvious from the context) we follow Mendelson's terminology and definitions of the "expressibility" ([Me64], §2, p117), and "representability" ([Me64], §2, p118), of number-theoretic relations and functions, respectively, in a

Now, a constructive view of Gödel's reasoning, as suggested in these essays, is that such a broad, and implicit (hence ambiguous), commitment is unnecessary, even if it is not formally invalid. It should, then, follow that, by the principle of Occams razor, we *ought not* to unrestrictedly assert [*R(x)*] as collectively satisfied by all *x* in the standard interpretation M of PA, although we *may* assert that [*R(x)*] is satisfied individually by any given *x* in M.

The significance of this, last, distinction is that the totality of values for which [*R(x)*] is satisfied in M is not, then, constructively definable as a formal mathematical object. In other words, the classical assumption that such values define a "set" in any Axiomatic Set Theory such as ZFC may, when interpreted constructively, introduce an anomaly, if not an inconsistency. Thus, we may not be able to make any constructively meaningful assertion about the totality of values that satisfy [R(x)] in M.

More precisely, we argue that every recursive[28] number-theoretic[29] relation need not be accepted as well-defining a (recursively enumerable[30]) sub-set of the natural numbers in any Axiomatic Set Theory that is a consistent[31] extension of PA[32], if we define a mathematical object and a set constructively as follows:

---

formal system such as P. However, following Gödel, we also refer to a number-theoretic relation as being "representable" in P when, strictly speaking, we mean that it is "expressible" in P as defined by Mendelson.

[28] We follow Mendelson's definitions of primitive recursive, and recursive functions and relations ([Me64], p120, 125). Primitive recursive functions / relations have the property that, for any set of natural number values given to their variables, the functions / relations are finitely computable / decidable on the basis of the axioms of any formal Peano Arithmetic.

[29] A number-theoretic function is one whose arguments and values are natural numbers, and a number-theoretic relation is a relation whose arguments are natural numbers ([Me64], p117).

[30] A recursively enumerable set is classically defined ([Me64], p250) as the range of some recursive number-theoretic function that, treated as a well-defined set, is implicitly assumed to be consistent with any Axiomatic Set Theory that is a model for PA.

[31] We define a language L as consistent if, and only if, there is no formula *F* of L such that both *F* and ~*F* (i.e., not *F*) are theorems in L, where we define *F* as a theorem of L if, and only if, there is a finite sequence

**Definition** (*i*): A *primitive mathematical object*, of a formal mathematical language, L, is any symbol for an individual constant, predicate letter, or a function letter (cf. [Me64], p46; also p1, p10), which is defined as a primitive symbol of L.

**Definition** (*ii*): A *formal mathematical object*, of a formal mathematical language, L, is any symbol for an individual constant, predicate letter, or a function letter that is either a primitive mathematical object of L, or one that can be introduced through definition (cf. [Me64], p82) into L without inviting inconsistency.[33]

**Definition** (*iii*): A *mathematical object*, of a formal mathematical language, L, is any symbol that is either a primitive, or a formal, mathematical object of L.

**Definition** (*iv*): A *set*, in the domain of any interpretation, M, of a formal mathematical language, L, is the range of any function in M that is an interpretation of a function letter that is a mathematical object of L.

---

of formulas of L, ending in $F$, such that each formula in the sequence is either an axiom of L, or a consequence of some of the formulas preceding it in the sequence by the deduction rules of the language.

[32] This means, loosely speaking, that we may not be able to give a set-theoretic definition of, say, Gödel's primitive recursive function $Sb(x\ v|Z(y))$, such that $\{x \mid x=Sb(y\ 19|Z(y))\}$ defines a set in any Axiomatic Set Theory with a Separation Axiom, or its equivalent, without introducing inconsistency.

Consequently, if we cannot consistently assume that every recursive function formally defines a recursively enumerable set, then it follows that we are unable to define a recursive set as a recursively enumerable set whose complement is also recursively enumerable. In some cases, there may be no such complement.

It further follows that, if $F(y)$ is an arithmetical function such that $F(k) = Sb(k\ 19|Z(k))$ for any given $k$, the assertion that the expression $\{x \mid x=F(y)\}$ defines a formal set by the Separation Axiom may require additional qualification.

[33] We note, incidentally, that, by defining a mathematical object precisely, the paradoxical element in the mathematical and logical "antinomies" is effectively eliminated; they define functions or relations that are not mathematical objects. Prima facie, except from a Platonic viewpoint, it thus seems of little significance whether such definitions are taken as defining concepts that are mathematically inconsistent (square circle), arguably inconsistent (Pegasus), logically inconsistent (Russell's impredicative set), or mathematical non-objects (the non-constructive elements of any axiomatic set theory).

Constructively, as a consequence of the arguments cited above, expressions such as "$(Ax)F(x)$", and "$(Ex)F(x)$", in an interpretation M of a formal theory P, may be taken to mean[34] "$F(x)$ is true for all $x$ in M", and "$F(x)$ is true for some $x$ in M", respectively, if, and only if, the predicate letter "$F$" is a formal mathematical object in P.

In the absence of a proof that "$F$" is such a mathematical object, the expressions "$(Ax)F(x)$" and "$(Ex)F(x)$" can, constructively, only be taken to mean that "$F(x)$ is true for any given $x$ in M", and "It is not true that $F(x)$ is false for any given $x$ in M", indicating that the predicate "$F(x)$" is well-defined, and decidable, for any given value of $x$, but that there may not be any uniformly effective method for such decidability.

## I-3    What exactly does this mean, and why is the distinction important?

The question arises: What exactly does this mean, and why is the distinction important?

Taking the latter part of the question first, the distinction is important because we do not, then, need to accept absolute limits on our ability to adequately formalise, and effectively communicate, mathematical concepts that are accepted as being common to a collective of gestalts[35].

Now, a central issue in the development of any significant Artificial Intelligence is that of understanding[36], and finding, effective methods of duplicating the cognitive and expressive processes of the human mind. This issue is being increasingly brought into sharper focus by the rapid advances in the experimental, behavioral, and computer sciences[37].

---

[34] We treat such, i.e., interpreted, symbolic expressions merely as abbreviations of semantically well-defined assertions of a language of communication in which we express our cognitive experiences.

[35] For instance, any recursively definable language would be considered as a mental concept that is common to a collective of gestalts.

[36] For instance, Manin [Ma02] draws attention to a paper in "Science", [DS99], that "summarizes some experimental results throwing light on the nature of mental representation of mathematical objects and physio-logical roots of divergences between, say, intuitionists and formalists", as below:

Penrose's "The Emperor's New Mind" and "Shadows of the Mind" highlight what is striking about the attempts, and struggles, of current work in these areas to express their observations adequately - necessarily in a predictable way[38] - within the standard interpretations of formal propositions as offered by classical theory.

So, the question arises: Are formal classical theories essentially unable to adequately express the extent and range of human cognition, or does the problem lie in the way formal theories are classically interpreted at the moment?

The former addresses the question of whether there are absolute limits on our capacity to express human cognition unambiguously; the latter, whether there are only temporal limits - not necessarily absolute - to the capacity of classical interpretations to communicate unambiguously that which we intended to capture within our formal expression.

Now, a thesis of these essays is that we may comfortably reject the first, by recognising that we can, indeed, constructively extend Tarski's definitions, of the "satisfiability" and

---

"[...] our results provide grounds for reconciling the divergent introspection of mathematicians by showing that even within the small domain of elementary arithmetic, multiple mental representations are used for different tasks. Exact arithmetic puts emphasis on language-specific representations and relies on a left inferior frontal circuit also used for generating associations between words. Symbolic arithmetic is a cultural invention specific to humans, and its development depended on the progressive improvement of number notation systems. [...]

Approximate arithmetic, in contrast, shows no dependence on language and relies primarily on a quantity representation implemented in visuo-spatial networks of the left and right parietal lobes."

[37] See, for instance ([RH01], footnote 18).

[38] We address this issue separately in §II-2, "Is mathematical truth verifiable effectively?".

"truth" of formal relations and propositions, respectively, under a given interpretation, verifiably.

## I-4    A constructive definition of Tarskian satisfiability

So, we return to the first part of the earlier question: What exactly does this mean?

Well, it means that we can, for instance, replace the strong[39], classical, implicit interpretation of Tarski's non-constructive definition, which is essentially the assertion:

> (*i*) **Classical satisfiability**: The well-formed formula [$F(x)$] of a formal system P is satisfied classically under an interpretation M of P if, and only if, the interpreted predicate $F(x)$ is satisfied collectively, *and* individually, by the elements in the domain of M,

by the weaker, explicit assertion:

> The well-formed formula [$F(x)$] of a formal system P is satisfied constructively under an interpretation M of P if, and only if, the interpreted predicate $F(x)$ is satisfied *either* collectively, *or* individually, by the elements in the domain of M.

We can, further, eliminate any implicit commitment to a Platonic ontology by making this definition constructive - in the sense of being effectively verifiable, as follows:

> (*ii*) **Individual satisfiability**: The well-formed formula [$F(x)$] of a formal system P is *individually* satisfied under an interpretation M of P if, and only if, given any value $k$ in M, there is an individually effective method (which may depend on the value $k$) to determine that the interpreted proposition $F(k)$ is satisfied in M.

---

[39] As remarked earlier, such a strong assertion may also be invalid under a constructive interpretation of Gödel's reasoning in Theorem VI of his 1931 paper [Go31a].

(*iii*) **Uniform satisfiability**: The well-formed formula [$F(x)$] of a formal system P is *uniformly*[40] satisfied under an interpretation M of P if, and only if, there is a uniformly effective method (necessarily independent of $x$) such that, given any value $k$ in M, it can determine that the interpreted proposition $F(k)$ is satisfied in M.

(*iv*) **Constructive satisfiability**: The well-formed formula [$F(x)$] of a formal system P is constructively satisfied under an interpretation M of P if, and only if, it is either uniformly satisfied in M, or it is individually satisfied in M.

Clearly, if [$F(x)$] is uniformly satisfied in M, then it is, obviously, individually satisfied in M. However, does the converse hold? More to the point, could Gödel's undecidable proposition, [$(Ax)R(x)$], be an instance of a formula [$R(x)$] that is individually satisfied in M (since Gödel shows that [$R(k)$] is provable in P for any numeral [$k$]), but not uniformly satisfied in M?[41]

An interesting consequence of an affirmative answer to this question would be that the interpreted arithmetical predicate $R(x)$ - which is instantiationally equivalent[42] to a primitive recursive relation - becomes Turing-undecidable[43]! Prima facie, this would appear to conflict with the classical postulation that a number-theoretic function is

---

[40] From here on, we prefer the term "uniformly" to "collectively", since the former is a familiar concept in various areas of mathematical reasoning, and the distinction - between individual satisfiability / truth and uniform satisfiability / truth - is very similar to the important distinction between, for instance, continuous functions ([Ru53], p65, §4.5) and uniformly continuous functions ([Ru53], p67, §4.13), and between convergent sequences ([Ru53], p115, §7.1) and uniformly convergent sequences ([Ru53], p118, §7.7).

[41] We address this issue separately in §II-3, "What is the significance of Gödel's First Incompleteness Theorem?".

[42] Two number-theoretic relations, say $f(x)$ and $g(x)$, are defined as instantiationally equivalent if, and only if, for any given natural number $n$, $f(n)$ holds if, and only if, $g(n)$ holds.

[43] We presume a degree of familiarity with the computability concepts introduced by Turing in his seminal paper [Tu36].

Turing-computable if, and only if, it is partial recursive[44]. However, the conflict may be illusory: the proof of equivalence seems to presume[45] that there is always a uniformly effective method for computing any given partial recursive function[46]. This proof would be invalid if we held, constructively, that Gödel has shown there are total[47] arithmetical functions, and relations, that may not be computable, or decidable, respectively, by any uniformly effective method[48].

## I-5     Extending Church's Thesis and defining effective computability

A key question, of course, is: why change?

One reason is that of verifiability. We are now able to define effective computability constructively, i.e., in a verifiable manner, by addressing the questions:

    (*v*) **Individual effective method**: When may we constructively assume that, given any sequence *s* of an interpretation M of P, there is an individually effective method to determine that *s* satisfies a given P-formula in M?

---

[44] A standard proof of this is given by Mendelson ([Me64], p233, Corollary 5.13 and p237, Corollary 5.15).

[45] This may, indeed, have been the presumption that allowed Turing [Tu36] to assert the essential equivalence between Gödel's definition of an undecidable, but true, arithmetical proposition, and his own definition of a number-theoretic Turing-uncomputable Halting function.

[46] See ([Me64], p226, Corollary 5.11).

[47] We define a number-theoretic function, or relation, as total if, and only if, it is individually computable, or individually decidable, respectively.

[48] In other words, although, for any given set of numerical values for its variables, a recursive function / relation is individually computable / decidable from the axioms of any Peano Arithmetic, there may be no common proof sequence within the Arithmetic that can compute / decide the function / relation for an arbitrary assignment of numerical values for its variables.

(*vi*) **Uniformly effective method**: When may we constructively assume that there is a uniformly effective method such that, given any sequence *s* of an interpretation M, *s* satisfies a given P-formula in M?

If the domain D of M can be assumed representable in P, then (*v*) can, indeed, be answered constructively; we simply extend the classical Church Thesis[49] as follows:

(*vii*) **Individual Church Thesis**: If, for a given relation $R(x)$, and any value *a* in some interpretation M of P, there is an individually effective method such that it will determine whether $R(a)$ holds in M or not, then every element of the domain D of M is the interpretation of some term of P, and there is some P-formula $[R'(x)]$ such that:

$R(a)$ holds in M if, and only if, $[R'(a)]$ is P-provable.

(In other words, the Individual Church Thesis postulates that the domain of a relation *R* that is effectively decidable individually in an interpretation M of some formal system P can only consist of mathematical objects, even if *R* is not, itself, a mathematical object.)

Moreover, (*vi*), too, can be answered constructively for any interpretation M of P, if we postulate:

(*viii*) **Uniform Church Thesis**: If, in some interpretation M of P, there is a uniformly effective method such that, for a given relation $R(x)$, and any value *a* in M, it will determine whether $R(a)$ holds in M or not, then $R(x)$ is the interpretation in M of a P-formula $[R(x)]$, and:

---

[49] We take the classical Church Thesis as the assertion that a number-theoretic function is effectively computable (partially) if, and only if, it is (partially) recursive ([Me64], p147, p227).

$R(a)$ holds in M if, and only if, $[R(a)]$ is P-provable.

(Thus, the Uniform Church Thesis postulates, first, that the domain of a relation $R$ that is effectively decidable uniformly in an interpretation M of a formal system P can only consist of mathematical objects; and, second, that $R$, too, is necessarily a mathematical object.)

## I-6    Some consequences, and a definition of *uniform* completeness

Some interesting consequences of the above are that:

(*ix*) The Uniform Church Thesis implies that a formula $[R]$ is constructively P-provable if, and only if, $[R]$ is uniformly satisfied in some interpretation M of P.

(*x*) The Uniform Church Thesis implies that if a number-theoretic relation $R(x)$ is uniformly satisfied in some interpretation M of P, then the predicate letter "$R$" is a formal mathematical object in P (i.e. it can be introduced through definition into P without inviting inconsistency).

(*xi*) The Uniform Church Thesis implies that, if a P-formula $[R]$ is uniformly satisfied in some interpretation M of P, then $[R]$ is uniformly satisfied in every interpretation of P.

(*xii*) The Uniform Church Thesis implies that if a formula $[R]$ is not P-provable, but $[R]$ is classically true under the standard interpretation, then $[R]$ is individually satisfied, but not uniformly satisfied, in the standard interpretation of P.

(*xiii*) The Uniform Church Thesis implies that Gödel's undecidable sentence GUS is individually satisfied, but not uniformly satisfied, in the standard interpretation of P.

(*xiv*) Under constructive extensions of Tarski's definitions of the satisfiability of the propositions of a formal mathematical language, P, under an interpretation, and of Church's Thesis, the domains of any two interpretations of P are denumerable (since they are representable in P) and isomorphic[50].

(*xv*) Under constructive extensions of Tarski's definitions of the satisfiability of the propositions of a formal mathematical language, P, under an interpretation, and of Church's Thesis, P is *uniformly* complete, in the sense that:

> **Definition** (*v*): A formal system P is *uniformly* complete if, and only if, a formula [*R*] is P-provable if, and only if, it is uniformly satisfied in some model of P.

We can also define effective computability, both individually and uniformly, along similar lines:

(*xvi*) **Individual computability**: A number-theoretic function $F(x)$ is individually computable if, and only if, given any natural number $k$, there is an individually effective method (which may depend on the value $k$) to compute $F(k)$.

---

[50] By "isomorphism", we mean the concept described by Mendelson [Me64], p90) as follows:

"We shall say that an interpretation M of the wfs of some first-order theory K is isomorphic with another interpretation M' of K if and only if there is a 1-1 correspondence g (called an isomorphism) of the domain D of M with the domain D' of M' such that:

(*i*) If $(A_j^n)^*$ and $(A_j^n)'$ are the interpretations in M and M', respectively, of $A_j^n$, then, for any $b_1, ..., b_n$ in D, $((A_j^n)^*(b_1, ..., b_n)$ if and only if $(A_j^n)'(g(b_1), ..., g(b_n))$.

(*ii*) If $(f_j^n)^*$ and $(f_j^n)'$ are the interpretations of $f_j^n$ in M and M', respectively, then, for any $b_1, ..., b_n$ in D, $g((f_j^n)^*(b_1, ..., b_n)) = (f_j^n)'(g(b_1), ..., g(b_n))$.

(*iii*) If $a_j^*$ and $a_j'$ are the interpretations of the individual constants $a_j$ in M and M', respectively, then $g(a_j)^* = a_j'$.

Notice that if M and M' are isomorphic, their domains must be of the same cardinality."

(*xvii*) **Uniform computability**: A number-theoretic function $F(x)$ is uniformly computable if, and only if, there is a uniformly effective method (necessarily independent of $x$) such that, given any natural number $k$, it can compute $F(k)$.

(*xviii*) **Effective computability**: A number-theoretic function is effectively computable if, and only if, it is either individually computable, or it is uniformly computable.

We can, then, give a constructive definition of uncomputable number-theoretic functions:

(*xix*) **Classical uncomputability**: A number-theoretic function $F(x_1, ..., x_n)$ in the standard interpretation M of P is classically uncomputable if, and only if, it is effectively computable individually, but not effectively computable uniformly.

This, last, removes the mysticism behind the fact that we can define a total number-theoretic Halting function that is, paradoxically, Turing-uncomputable. Ipso facto, it also validates some of the constructive objections (cf. [Hw98]) to Cantor's argument, that the cardinality[51] of the set of real numbers exceeds that of the set of integers.

If we accept the classical definition of a number as real if, and only if, it is the limit of a Cauchy sequence[52] of rational numbers[53], and that of a sequence of rationals as a function

---

[51] The cardinal number $C(Y)$ of a set $Y$ is defined to be the set of all sets $X$ that are equinumerous with $Y$ (i.e., for which there is a one-one correspondence between $Y$ and $X$). (cf. [Me64], p2).

[52] A numerical sequence $\{s_n\}$ is said to be a Cauchy sequence if, and only if, for every $e > 0$ there is an integer $N$ such that $n >= N$, $m >= N$ implies $|s_n - s_m| =< e$. (cf. [Ru53], p39, Definition, §3.10).

The importance of Cauchy sequences in mathematics is that a numerical sequence $\{s_n\}$ converges if, and only if, it is a Cauchy sequence. (cf. [Ru53], p39, Theorem, §3.11).

[53] This follows from the equivalence between the theorem that a sequence $\{s_n\}$ converges if, and only if, it is a Cauchy sequence ([Ru53], p39, §3.11), and the theorem that every Dedekind cut ([Ru53], p3, §1.4) defines a unique real number ([Ru53], p9, §1.32).

whose domain is the integers, and whose range is the rationals (i.e., the function defining the sequence is a mapping from the integers into the rationals), then Cantor's diagonal argument[54] simply shows that some real numbers are generated by classically uncomputable functions that are not effectively computable uniformly.

However, since, constructively, the function must be effectively computable individually, and there can only be denumerable effectively computable functions (by §(*xiv*) above), it follows that:

> (*xix*)(*a*) Under constructive extensions of Tarski's definitions of the satisfiability of the propositions of a formal mathematical language, P, under an interpretation, and of Church's Thesis, the real numbers are denumerable (hence Cantor's Theorem[55] does not hold).

---

[54] Cantor's diagonal argument is that if w e assume there is some sequence $\{f_n\}$ that maps the integers onto the real numbers (i.e., every integer occurs somewhere in the mapping) then we can define a real number that is not in the mapping - a contradiction. He concludes, therefore, that there can be no such mapping and, by definition, the cardinality of the reals must exceed that of the integers (in other words, in some sense, there must be more real numbers than there are integers).

However, constructively, the argument is unacceptable (cf. also [An03b]) since, if the sequence $\{f_n\}$ is effectively uncomputable uniformly, we cannot define a real number that is not in the mapping. In other words, since, given any $n$, we may have that $f_n$ is effectively computable individually, but there is no uniformly effective method such that, given any $n$, it will compute $f_n$, we cannot define a real number by Cantor's diagonal argument, as this, last, presumes that $\{f_n\}$ can be treated as uniformly computable.

[55] Cantor's (Power Set) Theorem ([Me64], Proposition 4.23, p183) is that, in any axiomatic set theory such as, say, ZFC, there can be no 1-1 correspondence between the members of the power set $P(x)$ of any set $x$ (intuitively, the set of all sub-sets of $x$), and the members of $x$. This theorem, however, appeals critically to the Separation Axiom of ZFC. Such an axiom admits the range of, both, uniformly computable, and individually computable, number-theoretic functions as mathematical objects of ZFC.

Classical theory simply terms the range of number-theoretic functions that are effectively computable uniformly as recursive sets that are also recursively enumerable, and the range of functions that are effectively computable individually, but not uniformly, as sets that are recursively enumerable, but not recursive.

Constructively, whilst the former sets may be validly assumed to be mathematical objects of ZFC, the latter may not, and any such assumption may invite inconsistency. This is the gist of the argument that we attempt to formalise in Meta-hypothesis 1 (the arguments for this are still tentative, and under revision) and its corollaries in §4 and §5 of Anand [An02c].

## I-7    The Uniform Church Thesis and the classical Church-Turing Theses[56]

The significance, of defining the satisfiability of a formula of P under interpretation explicitly in terms of individually, and uniformly, effective methods, and of expressing Church's Thesis constructively, is seen if we note that any computer can be designed to recognise a "looping" situation; it simply records[57] every instantaneous tape description at the execution of each machine instruction, and compares the current instantaneous tape description[58] with the record[59].

Now, we could instruct such a machine to assign arbitrary values to those undefined instances of the Halting predicate whose occurrences cause the machine to loop. Prima facie, this would define an individually effective decision method that is dependent entirely on the particular Halting function that is being computed, and cannot be predicted.

---

[56] This section presumes a degree of familiarity with the computability concepts introduced by Turing in his seminal paper [Tu36], and with classical interpretations of his reasoning as contained in, say, Mendelson ([Me64], p229-257).

[57] We may picture this operation as that of, say, a virtual teleprinter, with an infinite two-dimensional memory, that permanently records the instantaneous tape description of a Turing machine, after the execution of each machine instruction, as a new, finite, line array - and with facilities for comparing the current array with the previous record and interrupting the operation of the Turing machine at the onset of a loop (repeated array).

[58] "An instantaneous tape description describes the condition of the machine and the tape at a given moment. When read from left to right, the tape symbols in the description represent the symbols on the tape at the moment. The internal state $q_s$ in the description is the internal state of the machine at the moment, and the tape symbol occurring immediately to the right of $q_s$ in the tape description represents the symbol being scanned by the machine at the moment." ([Me64], p230, footnote 1).

[59] In principle, there is nothing in the definition of a Turing machine that bars such an operation, which would be that of a deterministic Turing-oracle, as discussed in §II-4.

Is such a machine a Turing machine? The classical answer to this question is not obvious if we do not appeal to the Turing Thesis[60].

However:

(*xx*) If we assume a Uniform Church Thesis, then every partial recursive number-theoretic function $F(x_1, ..., x_n)$ may have a unique constructive extension as a total function[61].

(*xxi*) If we assume a Uniform Church Thesis, then not every effectively computable function is classically Turing computable (so Turing's Thesis does not, then, hold).

(*xxii*) If we assume a Uniform Church Thesis, then not every (partially) recursive function is classically Turing-computable.

## II    Some significant issues in interpreting classical theory constructively

Although detailed proofs of the above consequences lie outside the scope of this essay[62], we review below, in a broader perspective, some of the significant issues, highlighted above, which indicate the need for balancing classical interpretations of Cantor's, Gödel's, Tarski's, and Turing's reasoning by an alternative, constructive and intuitionistically unobjectionable, interpretation of classical foundational concepts in which non-algorithmic truth is defined effectively.

---

[60] Turing's Thesis is that every effectively computable function is classically Turing computable ([Me64], p237). Although also referred to as the Church-Turing Thesis, CT, the latter is significantly different in that it implies that the Turing Thesis is equivalent to Church's Thesis. We consider the universal validity of Turing's Thesis, and of the equivalence, in §II-4, "What is the significance of Turing's Halting Theorem?".

[61] In other words, because the Turing machine that is computing the function loops for a given input, we cannot conclude that the function is undefined for that input. The value of the function, for the particular input, may simply require a method of computation that is unique to the input.

[62] Formal arguments, including proofs, of the above consequences of a constructive interpretation of classical theory, are expressed in Anand [An02c].

In particular, we address the questions:

(1) Are Platonism and Formalism incompatible doctrines?

(2) Is mathematical truth verifiable effectively?

(3) What is the significance of Gödel's First Incompleteness Theorem?

(4) What is the significance of Turing's Halting Theorem?

(5) Can all mental concepts be expressed mathematically?

(6) Can a constructive interpretation of Peano Arithmetic model some of the more paradoxical concepts of Quantum Mechanics?

## II-1    Are Platonism and Formalism incompatible doctrines?

Broadly reflecting the classical viewpoint:

**Platonism**[63]: *Platonism* may be viewed as the belief that abstract mathematical concepts are objective mathematical realities that can be perceived in a manner analogous to sense perception. Known also as mathematical realism, such a philosophy holds that mathematical entities exist independently of the human mind. Thus humans do not invent mathematics, but rather discover it, and any other intelligent beings in the universe would presumably do the same.

---

[63] These comments are based on the evolving, inter-active (yet, prima facie, astonishingly representative of standard classical theory), views on the Philosophy of Mathematics expressed in the open-source, on-line, encyclopedia Wikipedia [WiPM].

**Formalism**: *Formalism* may be viewed as the belief that mathematical statements may be thought of as statements about the consequences of certain string manipulation rules in the investigation of formal axiom systems. A theorem is not an absolute truth, but a relative one: *if* you assign meaning to the strings in such a way that the rules of the game become true (ie, true statements are assigned to the axioms and the rules of inference are truth-preserving), *then* you have to accept the theorem, or, rather, the interpretation you have given it must be a true statement.

That the two doctrines are, generally, seen as incompatible is reflected in Gower's remarks [Gw02]:

The basic Platonist position is rather simple. Mathematical concepts have an objective existence independent of us, and a statement such as "2+2=4" is true because two plus two really does equal four. In other words, for a Platonist mathematical statements are pretty similar to statements such as "that cup is on the table" even if mathematical objects are less tangible than physical ones. ...

Formalism is more or less the antithesis of Platonism. One can caricature it by saying that the formalist believes that mathematics is nothing but a few rules for replacing one system of meaningless symbols with another. If we start by writing down some axioms and deduce from them a theorem, then what we have done is correctly apply our replacement rules to the strings of symbols that represent the axioms and ended up with a string of symbols that represents the theorem. At the end of this process, what we know is not that the theorem is "true" or that some actually existing mathematical objects have a property of which we were previously unaware, but merely that a certain statement can be obtained from certain other statements by means of certain processes of manipulation.

That such apparent incompatibility may also have disconcerting, schizophrenic, undertones is reflected in the remark of a distinguished mathematician who, reportedly, described himself as a closet Platonist. On weekdays, when he was teaching (communicating to students), he lived in a world of formal structures; on weekends, when he was theorising (communing with himself), he lived in a world of Platonic ideals.

Unconsciously echoing this sentiment, set-theorist Saharon Shelah, too, reflected, in a soliloquy [Sh91]:

> "Does this mean you are a formalist in spite of earlier indications that you are Platonist?" I am in my heart a card-carrying Platonist seeing before my eyes the universe of sets, but I cannot discard the independence phenomena.

However, it can be argued, reasonably, that the philosophical paraphernalia carried by the two isms, Platonism and Formalism, may be largely irrelevant. The two may, simply, be complementary modes of thinking that we adopt, as convenient, to achieve our immediate purpose - the former as an aid to expression, the latter as an aid to communication.

Moreover, in practice, an implicit acceptance that they are two sides of the same coin may not be uncommon, despite debates on the subject that sometimes seem to suggest otherwise. In other words, every mathematician could be considered both as a creative Platonist, and an academic Formalist.

The question arises: Why, then, is there an unresolved debate on the issue?

Let us digress a little into some speculation. The ability to externally express, in some form, its internal states, can be seen as the indication, if not the definition, of a life form. Similarly, the ability to communicate such expression, to other life forms, can be seen as the indication, if not the definition, of intelligent life forms.

On this view, natural selection would simply be the reward of life forms whose intelligence is synonymous with the effectiveness of their ability to communicate. Under such a thesis (i.e., natural selection), such life forms could, possibly, have a greater chance of adapting holistically; and of retaining, for a longer period, a comparatively larger measure of a dynamically evolving identity, in an environment that is inimical to the survival of any concept of a static identity.

Viewed thus, the need, and ability, to express internal states externally could have a primordial significance that is rooted in the very ability of an evolutionary life form to survive, on the basis of instinctive reactions. In contrast, the need, and ability, to communicate effectively, would only be of significance consequently, in order to continuously evolve, and inter-act more harmoniously with its environment, through conscious action.

Returning to mathematics, the parallel would be that appreciation of the significance, and need, for Platonic modes of thought - as essential to the expression of abstract mathematical concepts, which can be viewed as existing within individual gestalts - could be at a more advanced stage, than appreciation of the significance, and need, for Formalist modes of thought. The latter may be significant only in the subsequent communication of such expression - for the identification of, say, that which can be accepted as uniformly common to a collective of individual gestalts.

So, one could, perhaps, say that Platonists, such as Gödel[64] earlier, and Shelah and Penrose today, are simply arbitrating, with reference to their individual gestalts, which of "Moreover, exactly as in the natural sciences ... inductio per enumerationem simplicem is by no means the only inductive method conceivable in mathematics. I admit that every mathematician has an inborn abhorrence to giving more than heuristic significance to such inductive arguments. I think, however, that this is due to the very prejudice that mathematical objects somehow have no real existence. If mathematics our mathematical propositions, when formally expressed in a precisely defined mathematical language, are intuitively true, and which are intuitively false, under their individual interpretations.

---

[64] Gödel's Platonism was explicit, and is well-documented by him. For instance, he remarks [Go51]:

By definition, such a view is, of course, unassailable. It may, moreover, be of value in understanding the cognitive processes of an individual brain, and of an individual mind.

However, clearly, such a view addresses only one half of human intellectual endeavor. This half would, first, be the attempt to individually express the state of the static synaptic elements, of the dynamically-evolving neuronic activity that is taken to represent an individual brain, within a symbolic language. And, second, it would be the subsequent attempt, to individually interpret, and relate, the symbols of a language to an existing state of static synaptic elements, of the dynamically-evolving neuronic activity that is taken to represent the individual brain.

The other half of human intellectual activity would, then, be that of determining which, of the concepts that are represented by such expressions, can be communicated uniformly in an unambiguous, and effective, manner that is independent of individual interpretations.

Curiously, although this area of intellectual activity would seem, prima facie, to be the natural preserve of Formalism, this issue does not seem to have been addressed, classically, as vigorously as it deserves, for reasons that seem to have a lot to do with the influence of Gödel's interpretation of his own formal reasoning.

---

describes an objective world just like physics, there is no reason why inductive methods should not be applied in mathematics just the same as in physics. The fact is that in mathemati cs we still have the same attitude today that in former times one had toward all science, namely, we try to derive everything by cogent proofs from the definitions (that is, in ontological terminology, from the essences of things). Perhaps this method, if it claims monopoly, is as wrong in mathematics as it was in physics."

For instance, most post-Gödelian discussions, on the foundations of mathematics, predominantly concern themselves with the question of whether we can, and how best we may, capture abstract mathematical concepts within formal languages that, it is implicitly accepted, cannot be interpreted uniquely.

Now, if we agree that abstract mathematical concepts may appear differently in different gestalts[65], then human experience in the experimental, observational and engineering disciplines indicates that we can, indeed, communicate the formalisation of mathematical concepts in an unambiguous, and effective, manner.

Prima facie, there is, thus, discordance between the theory and usage of mathematical languages. One would, reasonably, expect this to be an over-riding foundational concern of Formalism. This, however, does not seem to be the case at the moment.

To the contrary, a significant limitation of classical interpretations, of the formal reasoning and conclusions of classical first order theory - based primarily on the work of not only Gödel, but also that of Cantor, Tarski, and Turing - is an implicit, even if sometimes uncomfortable, acceptance of the argument - based, again, on an implicitly Platonic interpretation of Tarski's Theorem - that the truth of some propositions of a formal arithmetic, such as Peano Arithmetic, PA, under an interpretation M, is both non-algorithmic and essentially unverifiable constructively.

In other words, echoing Gödel's beliefs, it is accepted classically that PA cannot be taken as a faithful, and complete, formalisation of our intuitive, Dedekind, arithmetic; so, either (as standard interpretations of Gödel's reasoning and conclusions implicitly imply) such arithmetic is not formalisable in principle, or there is some, yet undiscovered, arithmetic that formalises M more faithfully[66].

---

[65] For instance, how would different individuals conceptualise a square circle?

However, such a view tends to gloss over a philosophically disturbing issue. The former is, intuitively, an unappealing, and implicitly self-limiting, admission; the latter, an unacceptable reflection on the competence of mathematicians to adequately select an appropriate set of primitive, axiomatic, assertions for PA as may be needed for PA to be an effective, and unambiguous, mathematical language of precise communication.

Now - since mathematics does serve increasingly as a language of precise expression and unambiguous communication for scientific discourse - a more reasonable thesis would be that such interpretations could, indeed, be balanced by an alternative, constructive and intuitionistically unobjectionable, interpretation of classical foundational concepts in which non-algorithmic truth is defined effectively[67].

---

[66] For instance, Gödel remarks [Go51]:

"However, as to subjective mathematics, it is not precluded that there should exist a finite rule producing all its evident axioms. However, if such a rule exists, ... we could never know with mathematical certainty that all propositions it produces are correct ...the assertion ... that they are all true could at most be known with empirical certainty ... there would exist absolutely unsolvable diophantine problems ..., where the epithet 'absolutely' means that they would be undecidable, not just within some particular axiomatic system, but by any mathematical proof the human mind can conceive."

[67] In other words, we argue that not every effective method is necessarily algorithmic, although every algorithm is an effective method. The possibility that mathematical tru th may be non-algorithmic, and yet constructive, is implicit in Gödel [Go51]:

"I wish to point out that one may conjecture the truth of a universal proposition (for example, that I shall be able to verify a certain property for any integer given to me) and at the same time conjecture that no general proof for this fact exists. It is easy to imagine situations in which both these conjectures would be very well founded. For the first half of it, this would, for example, be the case if the proposition in question were some equation $F(n) = G(n)$ of two number-theoretical functions which could be verified up to very great numbers $n$".

The possibility is also implicit in Turing's remarks ([ Tu36], §9, para II).

More precisely, it can be argued that some foundational concepts - implicitly accepted as intuitively unexceptionable in the classical interpretations of Cantor's, Gödel's, Tarski's and Turing's reasoning - can, under such thesis, be explicated effectively in non-Platonic interpretations that consider whether, and, if so, when and how, we may, within classical logic and without inviting inconsistency:

* define a mathematical object formally;

* define mathematical truth effectively:

* define effective methods of numerical computation non-algorithmically;

* differentiate between the concepts "Algorithmic truth: A formula $F(x)$ of a language L is uniformly true under an interpretation M for all $x$ in the domain of M" and "Non-algorithmic truth: A formula $F(x)$ of a language L is individually true under an interpretation M for any given $x$ in the domain of M";

* assert that two formulas of a formal system, under a given interpretation, have "the same meaning".

In other words, there may, arguably, be a need to free Formalism from the influence of an extended, post- Gödelian, Platonist philosophy - such as Realism. A need that arises from, and reflects, the failure of classical theory to distinguish between the representation, within a language, of those abstract concepts (mental constructs of an individual mind) that are individually significant within an individual gestalt, and th ose of these abstract concepts that are, further, communicable in an unambiguous and effective manner, and which may, therefore, be termed as uniformly significant within a collective of gestalts.

Prima facie, most of the challenges faced in unambiguously, and effectively, communicating mathematical concepts, seem to involve conclusions arrived at from debatable Realist premises[68]. Premises that, when introduced authoritatively into formal mathematical reasoning, and into scientific discourse, permit us to logically validate our subjective intuitive perceptions as being reflective of some absolute Truth that - contrary to our intuitive experience - must be of universal significance in a Utopian, Platonic world.

We could even go further, and consider whether the relation of a particular language (whether mathematical or not), to that which the language seeks to express, should be the preserve of the philosophers of the language, and of those who study the nature of the mind and of consciousness - not of the logicians or mathematicians who simply provide the tools for such study.

In other words, the selection of the alphabet, selection of rules for defining well-formed words, phrases, and declarative atomic and compound sentences, selection of primitive objects (constants) and primitive truths (axioms), selection of rules of deduction for assigning truth values to non-axioms, etc., in a language should be consequent to the adequate resolution of philosophical issues that consider the question of what we visualise as, and how we visualise, the properties and relations between elements of an abstract, Platonically conceived, ontology- in our individual gestalts - that a logician or mathematician is being asked to represent within a formal language.

Whether the above properties and relations between elements of a chosen abstract, Platonically conceived, ontology have limited, individual significance, or have a wider, common significance, should also be the rightful preserve of philosophy in general, along

---

[68] For instance, that mathematical truth is discovered.

with questions concerning the nature of such elements, their properties, and their relationships.

However, once we have defined a language, we should not be at liberty to use the rules of the language to deduce the existence of new, unique, elements of the ontology (as opposed to elements of the language) that are not explicitly specified by the definitions of the language. Such deduction would amount to a creation by extraneous definition that would contradict the premise that the language has already been specifically defined to express selected pre-conceived, Platonic, concepts.

In other words, the aim of mathematics should not be to introduce such extraneous definitions into a language, but to study the consequences of a mathematical language that, beyond the considerations involved in its definition, are independent of the Platonic concepts that the language was designed to express. By this yardstick, standard interpretations of Cantor's diagonal argument[69], or of his power set theorem (as also of some of the transfinite elements of set theory) seem to extravagantly violate this principle.

## II-2   Is mathematical truth verifiable effectively?

In a general talk, Gowers [Gw02] remarked that:

> If you ask a philosopher what the main problems are in the philosophy of mathematics, then the following two are likely to come up: what is the status of mathematical truth, and what is the nature of mathematical objects? That is, what gives mathematical statements their aura of infallibility, and what on earth are these statements about?

---

[69] We address some of the non-constructive issues involved in classical interpretations of Cantor's diagonal argument in [An03b].

Thus, even 70 years after Gödel's seminal 1931 paper [Go31a], where he highlighted that there is an essential asymmetry between classically true, and classically provable, arithmetical propositions, the roots of such distinction continue, apparently, to elude effective expression.

Now, prima facie, for any language to be termed as a language of unambiguous, and effective, communication, the truth of its propositions, under any interpretation, ought to be unambiguous, and effectively verifiable independently of the domain of the interpretation.

We are, thus, faced with the questions: Is mathematical truth necessarily unverifiable effectively, and are there theoretically absolute limitations on unambiguous, and effective, communication?

Now, following Gödel's interpretation of his own reasoning, and conclusions, in his 1931 paper [Go31a], standard interpretations of classical mathematical theory seem to implicitly imply that, even in a mathematical language as basic as formal Peano Arithmetic, the most fundamental of our[70] intuitive mathematical concepts cannot be expressed, and communicated in a complete, unambiguous, and effective manner[71].

---

[70] If we accept that an awareness of the concept of finite counting is also exhibited in the behaviour of other species, then the concept may be intuitively fundamental in a far larger sense.

[71] Which gives an unintended, and devilishly misleading, twist to Russell's (justly) famous quip [Rs01]:

"Thus mathematics may be defined as the subject in which we never know what we are talking about, nor whether what we are saying is true."

The extent to which such interpretations are accepted as setting implicit, and, apparently absolute, limitations on unambiguous, and effective, communication, depends, no doubt, on individual psychological, and philosophical, predilections. Thus, although professional mathematicians and mathematical logicians might succeed in placing such limitations in comfortably abstract, albeit counter-intuitive, perspective, such comfort may be denied to other disciplines that primarily need to express their fundamental concepts and observations in an unambiguous, and effectively communicable, manner. For instance, seekers of extra-terrestrial intelligence, or those striving to replicate human intelligence in controllable artifacts, could find such acceptance disturbingly constricting, and philosophically disquieting.

In this essay, we suggest that such limitations are not absolute, but self-imposed. We argue that they are rooted in a removable ambiguity in the standard interpretations of Tarski's, non-verifiably formulated, definitions of the satisfiability, and truth, of the formulas of a formal system, under a given interpretation (essentially, in an unspecified, and implicitly Platonic, domain).

Now, a question of increasing interest to scientific disciplines that look to mathematics for providing a language of reliable, and verifiable, external expression and communication[72], is: Are the concepts "non-algorithmic" and "non-constructive" necessarily synonymous in classical logic and mathematics?

Gödel argues that, in a formal language as basic as Peano Arithmetic, there are undecidable sentences that can be recognised as true under classical interpretation, but which are not provable within the system. Does this imply that such recognition, in some cases, cannot be duplicated in any artificially constructed and, more important, artificially controlled, mechanism or organism whose design is based on classical logic?[73]

---

[72]For most scientific disciplines, the authority of the standard interpretations of classical mathematics is seen, and accepted (perhaps with some element of reluctance, since such acceptance occasionally flies against the grain of observation and experience) not only as absolute, but also as implicitly promising sufficiency, when needed, to help bridge the seemingly unbridgeable chasm that sometimes confronts such disciplines - between a Platonic world of abstract objects, and the real world of sensory perceptions!

The scientific, and philosophical, dimensions of an affirmative answer to the last question have been broadly reviewed, and addressed, by Penrose in [Pe90] and [Pe94]. Penrose's argument is based on a strongly Platonist thesis that sensory perceptions simply mirror aspects of a universe that exists, and will continue to exist, independent of any observer ([Pe90], p123, p146)[74]. On this view, individual consciousness would be a discovery of what there is (cf. [Pe90], p124), and be independent of the language in which such discovery is expressed. It follows that recognition of intuitive truth would be individually asserted - and, implicitly, fallible - correlations between the unverifiable - and, ipso facto, infallible - intuitive experiences of an individual consciousness, and the formal expressions of a communicable language.

The issue, then, is whether classical logic can adequately formalise intuitive truth to make it infallible, or whether such recognition is essentially fallible[75]. Penrose opts for the inadequacy of classical logic to completely capture a Platonic mathematical reality that, he believes, manifests itself, first, in thought - which originates in the mind consequent to sensory experience - and, second, in its representation in a communicable language. He supports his view by highlighting the "ethereal" presence, and non-verifiable properties, of various non-algorithmic ([Pe90], p168), and implicitly non-constructive, mathematical concepts that are accepted in our formal languages as essential to classical mathematics ([Pe90], p123-8).

---

[73] We note that the question may have economic significance globally, particularly in areas relating to the development of strategic and infra-structural products, facilities and services that are based on the proposed replication of human intelligence by artificial mechanisms or organisms.

[74] An obvious, but arguably relevant, objection to this argument is that it assumes multiple, spatially separated, observers can each, Deity-like, acquire identical knowledge of a Universe simultaneously without altering, or even indirectly influencing, the knowledge that is sought to be acquired.

[75] Of course, there is an inescapable element of circularity in considering the fallibility of assertions that are asserted as intuitively true.

Although Penrose's arguments represent only one, and perhaps an arguably extreme, point of view[76], they emphasise that classical mathematics may yet need to adequately legitimise the acceptance, into a theory, of formally definable mathematical objects (cf. [Pe90], p147), most obviously those that can be argued as being essentially non-constructive.

Now we note that Penrose appears to base his thesis on, amongst others, a classical consequence of Gödel's reasoning and conclusions:

> We cannot express Tarskian definitions of the "satisfiability" and "truth" of formal expressions under an interpretation algorithmically ([Pe90], p159).

This is, essentially, an intuitive interpretation of:

> **Tarski's Theorem**: The set $Tr$ of Gödel-numbers of the formal expressions of a first order Peano Arithmetic that are true in the standard model is not arithmetical. In other words, there is no formula $[F(x)]$ in any formal Peano Arithmetic such that $Tr$ is the set of numbers $k$ for which $[F(k)]$ is true in the standard model.

Penrose concludes that, although we may follow a common intuitive process for discovering common mathematical aspects of the universe, not all our mathematically expressible discoveries are expressible by classical algorithms ([Pe90], p533, p548).

---

[76] See *Psyche*, Vol. 2(9), June 1995, *Symposium on Roger Penrose's Shadows of the Mind*.

However, Penrose's arguments also appear to imply further, albeit implicitly, that our recognition of intuitive "arithmetical truth" - even when this is accepted as being adequately formalised by the classical Tarskian definitions of the "satisfiability" and "truth" of formal expressions under an interpretation - is "absolutely" non-constructive (cf. [Pe90], p145-6).

Thus, although Penrose does not seem to question the mathematical form of Church's Thesis ([Pe90], p64-65) - which, essentially, postulates that every effectively computable function is algorithmic - he seems to conclude from his arguments, concerning the inadequacy of classical logic, that there may be non-algorithmic, non-constructive, ways of acquiring mathematical insight and knowledge ( [Pe90], p538)[77].

As is evidenced in his discussion of Lucas' Gödelian argument[78] [Lu61], Penrose does not appear ([Pe90], p539) to entertain the possibility that there may be non-algorithmic reasoning that could be intuitionistically accepted as constructive; his arguments seem to, implicitly, treat the terms "non-algorithmic" and "non-constructive" as synonymous[79].

---

[77] Penrose's belief, that non-algorithmic effective methods are necessarily non-constructive in classical mathematics, may, however, be reflective of the thinking amongst a wide set of scientific disciplines, as seen in Gurney's remarks [ Gu96], (reproduced in §II-5), in the context of neural nets.

[78] "Gödel's theorem seems to me to prove that Mechanism is false, that is, that minds cannot be explained as machines."

[79] Davis [Da90] critically reviews this particular aspect of Penrose's argument. He argues that: "... Gödel's incompleteness theorem (in a strengthened form based on work of J.B. Rosser as well as the solution of Hilbert's tenth problem) may be stated as follows: There is an algorithm which, given any consistent set of axioms, will output a polynomial equation P = 0 which in fact has no integer solutions, but such that this fact can not be deduced from the given axioms. Here then is the true but unprovable Gödel sentence on which Penrose relies and in a particularly simple form at that. Note that the sentence is provided by an algorithm. If insight is involved, it must be in convincing oneself that the given axioms are indeed consistent, since otherwise we will have no reason to believe that the Gödel sentence is true".

However, the real issue is not whether there is an algorithm that outputs P = 0, but whether, for any given set of natural number values for its free variables, the fact that P = 0 has no integer solutions can be determined in a non-algorithmic, yet constructive way.

However, what Penrose, for instance, views as the essentially, and absolutely non-constructive, aspects of mathematical concepts, may simply be manifestations of a *removable* ambiguity in the classical Tarskian definitions of the satisfiability, and truth, of our formal expressions under an interpretation.

In other words, we may simply need to:

(*i*) first, remove possible ambiguities in Tarski's definitions by explicitly introducing the concept of effective verification into Tarski's definitions, and,

(*ii*) second, explicitly specify when a method of verification can be considered effective by re-formulating Church's Thesis.

An interesting consequence, and raison d'etre, of such an alternative, constructive and intuitionistically unobjectionable, interpretation of Gödel's reasoning, and of Tarski's definitions of satisfiability and truth, is that Peano Arithmetic can, then, be shown to be complete, albeit in a broader sense than that in which the term is defined classically.

Thus, we could consider a language *verifiably complete* if, and only if, every true statement of the language is effectively verifiable under definitions such as, say, the following:

(*a*) **Individual truth**: A formula $[F(x)]$ of a formal system P is individually true under an interpretation M of P if, and only if, given any value $k$ in M, there is an individually effective method (which may depend on the value $k$) to determine that the interpreted proposition $F(k)$ is satisfied in M.

(*b*) **Uniform truth**: A formula $[F(x)]$ of a formal system P is uniformly true under an interpretation M of P if, and only if, there is a uniformly effective method (necessarily independent of *x*) such that, given any value *k* in M, it can determine that the interpreted proposition *F*(*k*) is satisfied in M.

(*c*)    **Effective truth**: A formula [*F*(*x*)] of a formal system P is effectively true under an interpretation M of P if, and only if, it is either individually true in M, or it is uniformly true in M.

We note that, if we eliminate the references to effective methods in the above, we arrive at Tarski's corresponding definitions of the truth of formal propositions under an interpretation (cf. [Me64], p49-52). However, clearly, these would no longer be able to distinguish between individual truth (satisfiability), and uniform truth (satisfiability), as suggested above.

Now, Gödel argues that, in any consistent, formal, system P that formulates Peano's Arithmetic, we can construct a valid expression of the system, say [*R*(*x*)], such that [*R*(*n*)] is P-provable for any given numeral [*n*], but [(A*x*)*R*(*x*)] is not P-provable. The classical interpretation of this is that although [(A*x*)*R*(*x*)] is not P-provable, it is true under its standard interpretation by Tarski's definitions of satisfiability and truth.

However, by implications that can be considered as implicit in Tarski's definitions, [*R*(*n*)] may be viewed alternatively as an expression whose standard interpretation, *R*(*n*), can be effectively asserted as holding individually - and not necessarily algorithmically - for any given natural number *n*, but *R*(*x*) cannot be effectively asserted as holding uniformly - in the sense of algorithmically - for all natural numbers *x* collectively.

In other words, the classical interpretations of Tarski's definitions of satisfiability and truth seem to contain an ambiguity: they implicitly imply the existence of an ambiguous effective method for deciding whether formal expressions such as [*R*(*x*)] are satisfied under a given interpretation. Specifically, they fail to entertain the possibility that such a method may be non-algorithmic.

Thus, for any given value $n$ of its free variable under a given interpretation, there may always be a - possibly $n$-specific - (hence non-algorithmic) method that can effectively decide whether the interpretation $R(n)$ of a formal expression such as $[R(n)]$ holds individually, even when there is no $n$-independent (algorithmic) effective method that can effectively decide whether the expression $[R(x)]$ is satisfied uniformly, under a given interpretation, when we substitute any numeral $[n]$ for its free variable.

The ambiguity in Tarski's definitions is reflected in reasoning such as:

> $[(Ex)F(x)]$ is provable in Peano Arithmetic. Hence, we can conclude the existence of some, unspecified, natural number $s$ such that $[F(s)]$ is provable in PA, and use $[F(s)]$ as a provable assertion in formal arguments.

Now, even if we accept the classical Tarskian definitions of the satisfiability and truth of number-theoretic assertions as definitive, a reasonable interpretation of "$[(Ax)F(x)]$ is provable in PA" could be an inductive one. In other words we could still interpret it as the assertion that there is some uniformly effective method, independent of $x$, of determining that, given any natural number $s$, $[F(s)]$ is provable in PA.

It follows that "$[(Ex)F(x)]$ is provable in PA", which is defined as "$[\sim(Ax)\sim F(x)]$ is provable in PA", would then be the assertion that there is no uniformly effective method, independent of $x$, of determining that, given any natural number $s$, $[\sim F(s)]$ is provable in PA.

However, for any given $s$, there may yet be some individual effective method, which depends on $s$, such that $[\sim F(s)]$ is provable in PA[80]. Can we, therefore, validly conclude - even classically - that $[F(s)]$ is necessarily provable in PA for some $s$?

Leaving the question unaddressed tacitly tolerates extravagant existential interpretations in classical theory. Thus, as Christer-Hennix notes [He04]:

> "...in 1959, A. Heyting characterized, not without chagrin, the development after 1931, at the Warsaw Conference, with the following words:
>
>> The danger to which I alluded just now, namely that a notion which was meant to be constructive, is afterwards interpreted in a non-constructive way, is not imaginary; on the contrary, several meta-mathematical notions have suffered by it. In general, those mathematicians who introduce the non-constructive interpretation, are not aware that they falsify the notion as it was originally intended.

Further, ignoring the issue also seems to obscure the raison d'etre of a mathematical language - to communicate unambiguously, and effectively. Thus, if mathematics is to serve both as a language that expresses all possible abstract (including Platonic) mathematical concepts in formal languages, and also to serve as a language of precise expression and unambiguous communication, then, particularly in the latter case, we may need to specify effective decision procedures for determining whether, or not, a proposition of a formal language of mathematical communication, L, is to be termed as true or not under each given interpretation M.

---

[80] This possibility, prima facie, admits (without rejecting the Law of the Excluded Middle) Brouwer's, and the intuitionist, objection to the general validity of putative instantiations that are concluded, classically, from interpretations of the existential theorems of a formal system.

For instance, similar to the familiar Church Thesis (that we may treat a number-theoretic function as effectively computable if, and only if, it is recursive) we may need to specify when, and only when, we can treat the truth of a number-theoretic relation as effectively decidable. The Individual and Uniform Church Theses, defined above in §I-5, could, indeed, serve as possible paths towards such an end.

## II-3  What is the significance of Gödel's First Incompleteness Theorem?

If we define satisfiability and truth effectively as suggested earlier, it turns out that, by Gödel's reasoning, every system of Peano Arithmetic is omega-inconsistent (see the definition below), and so Gödel's Theorem VI (the famous First Incompleteness Theorem of his seminal 1931 paper [Go31a]) holds vacuously.

To see this, consider Gödel's recursive relation $Q(x, p)$, which is represented in P by $[R(x, p)]$, where $[(Ax)R(x, p)]$ is Gödel's "undecidable" proposition. Thus, although we can effectively conclude, from Gödel's reasoning in the first part of his proof of Theorem VI, that "$Q(n, p)$ holds individually for any given natural number $n$", we cannot, prima facie, assume that this is equivalent to the non-constructive, infinite, compound, sentence "$Q(x, p)$ holds uniformly for all natural numbers $x$".

The distinction between the two assertions may be better expressed in terms of classical Turing machines. Thus, given any natural number $n$, it follows - from Gödel's reasoning that $[R(n, p)]$ is P-provable - that there is always some effective method that will terminate in a finite, even if indeterminate, number of steps $t(n)$ if, and only if, $R(n, p)$ holds.

However, since $[(Ax)R(x, p)]$ is not P-provable, there may not be any Turing machine such that, given any $n$, it will halt in a determinate number of steps, $t(n)$, if, and only if, $R(n, p)$ holds . In other words, there may be no classical Turing machine that computes the function $t(x)$; in Turing's terminology, $t(x)$ may be Turing-uncomputable, even though $t(n)$ is effectively computable for any given $n$. Thus, a constructive interpretation of Gödel's reasoning implies that the Turing thesis is false.

Now, the thesis of a constructive interpretation of Gödel's reasoning and conclusions is that we may not interpret, for instance, the meta-assertion "PA proves: $[(Ax)F(x)]$ " as the non-verifiable, Tarskian meta-assertion:

$F(x)$ is satisfied by any given element $x$ of the domain of M .

We must interpret it, instead, as the verifiable meta-assertion:

There is a uniformly effective method (algorithm/Turing machine) such that, given any element $x$ of the domain of M, it will effectively decide that $F(x)$ is satisfied in M.

It follows that both the meta-assertions, "PA does not prove $[(Ax)F(x)]$", and "PA proves $[\sim(Ax)F(x)]$", then interpret constructively as the meta-assertion:

There is no uniformly effective method (algorithm/Turing machine) such that, given any element $x$ of M, it will effectively decide that $F(x)$ is satisfied in M.

Consequently, a constructive interpretation of Gödel's reasoning and conclusions implies that there can be no undecidable propositions in PA; in other words, that PA is complete in the sense that the truth of any arithmetical proposition is effectively, i.e., either individually, or uniformly, verifiable in PA.

However, we are, then, immediately faced with the question: Since the standard interpretation of Gödel's reasoning and conclusions asserts that PA is classically incomplete, how definitive is the standard interpretation?

Now, in Theorem VI of his 1931 paper, Gödel essentially argues that his "undecidable" proposition, [(A*x*)*R*(*x*)], is such that:

If [(A*x*)*R*(*x*)] is PA-provable, then [~*R*(*n*)] is PA-provable for some numeral [*n*].

Now, by standard logical axioms, we have that:

If [~*R*(*n*)] is PA-provable for some numeral [*n*], then [~(A*x*)*R*(*x*)] is PA-provable.

It thus follows that Gödel has, essentially, argued that:

If [(A*x*)*R*(*x*)] is PA-provable, then [~(A*x*)*R*(*x*)] is PA-provable.

Clearly, it should now follow, by the standard Deduction Theorem of first order logic, that:

[(A*x*)*R*(*x*) => ~(A*x*)*R*(*x*)] is PA-provable,

and so:

[~(A*x*)*R*(*x*)] is PA-provable.

However, at this point, standard interpretations of Gödel's reasoning appeal to his explicit assumption that PA is omega-consistent[81] in order to conclude that the PA-provability of [~(A*x*)*R*(*x*)] cannot be inferred from the above meta-argument.

**Omega-consistency**: A formal system P is omega-consistent if, and only if, there is no P-provable formula [~*F*(*x*)] such that [*F*(*s*)] is P-provable for any given numeral [*s*].

---

[81] The concept of omega-consistency was introduced by Gödel (Go31a], p23-24). Although Gödel's reasoning was specific to his formal system P, he noted that it could be replicated equally within other formal systems of Set Theory, such as ZFC, and Peano Arithmetics, such as standard PA.

Now, unless the omega-consistency of PA has some deeper, intuitive significance philosophically, this is not a reasonable inference. Since the Deduction Theorem is a fundamental theorem of classical logic, we must, using Occam's razor, conclude from Gödel's reasoning, first, that $[\sim(Ax)R(x)]$ is PA-provable, and, second, that PA is omega-inconsistent - and so Gödel's Theorem VI holds vacuously[82]!

Now we note that the omega-inconsistency of PA is, actually, natural, and intuitively unobjectionable, under a constructive interpretation of the concept of "PA proves: $[(Ax)F(x)]$" as described earlier.

Under such interpretation, an omega-inconsistent PA does not imply that PA, or any of its interpretations, are either inconsistent, or unnaturally consistent; it simply implies that there are arithmetical relations that cannot be verified uniformly by any effective method (algorithm / Turing machine) over the domain of their interpretation.

The above suggests that it may be the absence of an adequately technical counter-argument that leaves Wittgenstein's viewpoint [Wi78] - and, possibly, that of others who have shared his reservations and misgivings, such as the ultra -intuitionist Yessenin-Volpin today (cf. [He04]) - vulnerable to the arguments advanced by classical interpretations of Gödel's reasoning and conclusion; the latter, implicitly, imply that any interpretation of Gödel's reasoning and conclusion must be accepted as essentially counter-intuitive on the basis of purely technical considerations.

Prima facie, the classical interpretation of Gödel's reasoning and conclusions seems strengthened by Rosser's argument that Gödelian undecidability can be established in a simply consistent PA without assuming omega-consistency. However, if a standard (textbook) exposition of Rosser's proof, such as in Mendelson [Me64], can be accepted as definitive, then Rosser's argument may be non-constructive, and intuitionistically objectionable[83].

---

[82] This argument is considered in detail in the web essays [An02b] and [An03i]. We note that some of Wittgenstein's remarks (cf. [An03h]) indicate that, prima facie, he, too, saw no intuitively significant philosophical grounds for treating the omega-consistency of Peano Arithmetic as a primitive concept, and for allowing it to over-ride an application of the Deduction Theorem.

---

[83] See [An02a] for a detailed analysis of Rosser's argument. Briefly, Rosser's informal meta-argument seems to be that, assuming [$H$] is a provable PA-formula, and that $j$ is the Gödel-number of a proof of [$H$] in PA, and given that:

(*i*) PA proves: [$H$] => PA proves: [$q(\underline{j})$],

where [$q(\underline{y})$] is also a PA-formula, we may conclude, by the Deduction Theorem, that:

(*ii*) PA proves: [$H$] => [$q(\underline{j})$].

However, such reasoning is invalid, since (*ii*) may not be a theorem of PA.

To see this, note that, given any natural number $h$, we could construct a Turing-machine that would decompose $h$ to check whether it is in fact the Gödel-number of some well-formed formula [$H$] of PA. Then, assuming the provability of [$H$], the program would conclude that the meta-assertion (*ii*) is invalid, since it does not hold for every natural number $j$.

In other words, the formal meta-mathematical expression of (*i*) is, actually, the meta-statement:

(*iii*) PA proves: [$H$] => (E$j$)PA proves: [$q(\underline{j})$].

Now, from this, we cannot conclude (as Rosser does) that:

(*iv*) PA proves: [$H$] => PA proves: [(E$j$)$q(\underline{j})$].

Thus, from (*iii*), we may only conclude that, under the standard interpretation, it is true there must be the number $j$. However, for $j$ to be introduced as a symbol into a formal PA-proof, as in Rosser's arguments (*i*) and (*ii*), where he applies the Deduction Theorem, we also need a formal PA-proof that the number $j$ does exist in every interpretation. However, [(E$j$)$q(\underline{j})$] may be PA-unprovable.

In other words, if $j$ is Turing-uncomputable, then $j$ cannot be constructively defined by any PA-formula, and so the Deduction Theorem cannot be applied to the formal meta-assertion (*iii*), to conclude (*ii*) within a formal proof sequence of PA, so as to yield Rosser's undecidable proposition RUS in PA. Of course, PA proves: [$q(\underline{j})$] would follow from standard logical axioms if it could be shown that [(E!$j$)$q(\underline{j})$] is PA-provable (where "!" signifies uniqueness).

Now Rosser's argument implicitly assumes that $j$ is always Turing-computable (hence there is always a PA-formula that uniquely represents $j$), and arrives at a contradiction. However, this only establishes that $j$ is Turing-uncomputable. In the absence of a PA-proof of (*iv*), we cannot logically conclude from this, as he seems to do, that there is no such $j$.

## II-4    What is the significance of Turing's Halting Theorem?

Consider the following form of the Halting Theorem:

> **Theorem 1**: There is no Turing machine *WillHalt* which, given an input string M+ *w*, will halt and accept the string if Turing machine M halts on input *w*, and will halt and reject the string if Turing machine M does not halt on input *w*, where, viewed as a Boolean function, *WillHalt*(M, *w*) (halts and) returns true in the first case, and (halts and) returns false in the second.

> **Proof**: Assuming that *WillHalt*(M, *w*) exists leads to an immediate contradiction.

However, the above does not address an implicit, and foundationally fundamental, non-constructive issue: What, precisely, do we mean by "M does not halt on input *w*"?

If *WillHalt*(M, *w*) could effectively determine somehow that M does not halt on input *w*, then this procedure could be built into M, and so M would halt on *w* (possibly with the perplexing annotation: "Sorry, this program is being halted as the program has determined that it does not halt!"). Thus, we cannot even define the machine *WillHalt*(M, *w*) as a Turing machine without inviting an immediate inconsistency [84]!

So, if we accept that *WillHalt*(M, *w*) is a well-defined Boolean function, then the significant conclusion from the above (and Turing's) argument is not that *WillHalt*(M, *w*) is Turing-uncomputable, since we cannot define a Turing machine *WillHalt*, but that, if the function is effectively computable, then the Turing Thesis is false.

---

[84] Prima facie, such definitions appear similar to that of a *Liar* proposition as "The *Liar* proposition is a lie", or that of a *Russell* set as {*x* | *x* is a member of the *Russell* set if and only *x* is not a member of *x*}!

Now, in his 1936 paper [Tu36], Turing did not conclude from his Halting argument that a function such as *WillHalt*(M, *w*) is not computable; he only concluded that such a function cannot be computed by a Logical Computing (Turing) Machine as defined by him. However, since he was of the view that effective computability was, indeed, equivalent to computability by an LCM, Turing later proposed that a function such as *WillHalt* could, perhaps, be defined as an "oracle" that is not necessarily either deterministic or mechanical.

As noted by Hodges [Ho00]:

> Turing's 1938 Princeton Ph.D. thesis, work conducted in close cooperation with Church, was entitled Systems of logic defined by ordinals, and published as (Turing 1939). Predominantly the work consisted of highly technical developments within mathematical logic. However the driving force lay in the question: what is the consequence of supplementing a formal system with uncomputable deductive steps? In pursuit of this question, Turing introduced the definition of an "oracle" which can supply on demand the answer to the halting problem for every Turing machine. ... Turing defined the "oracle" purely mathematically as an uncomputable function, and said, 'We shall not go any further into the nature of this oracle apart from saying that it cannot be a machine.' The essential point of the oracle is that it performs non-mechanical steps.

Clearly, if the Turing Thesis is false, then *WillHalt* can be a machine; further, if *WillHalt* is effectively computable, then Turing's Thesis is false. So, the point to consider is whether there are mechanical computational processes that are not obviously duplicatable by the operations of a classical Turing machine.

Now, we can design a mechanical computer that will recognise a "looping" situation; it simply records every instantaneous tape description at the execution of each machine instruction, and compares the current instantaneous tape description with the record (we may visualise this as the provision of a deterministic oracle, in the form of a virtual teleprinter, with an interrupt facility). We could then instruct the machine to assign arbitrary values to those undefined instances of *WillHalt*(M, *w*) whose occurrences cause the machine to loop.

Clearly, if we can show, or accept, that the extended *WillHalt*(M, *w*) is well-defined, then Turing's Thesis is false. For instance, consider the following:

**Theorem 2**: The Uniform Church Thesis implies that the Halting problem is effectively solvable.

**Proof**: If a number-theoretic function $F(x)$ is Turing-computable, then it is partial recursive ([Me64], p233, Corollary 5.13). We may thus assume that $F$ is obtained from a recursive function $G$ by means of the unrestricted *mu*-operator; in other words, that $F(x) = mu_y(G(x, y) = 0)$ ([Me64], p214).

If, now, [$H(x, y)$] expresses $\sim(G(x, y) = 0)$ in a formal system of Peano Arithmetic such as standard PA, we then consider the PA-provability, and truth in the standard interpretation M of PA, of the formula [$H(a, y)$] for a given numeral [$a$] of PA.

We consider the argument:

(*a*) Let $Q_1$ be the meta-assertion that [$H(a, y)$] is not effectively true in M. Hence there is no effective method in M to determine that, for any given *y* in M, *y* satisfies [$H(a, y)$] classically. It follows that there is no uniformly effective method (Turing machine) in M to determine that, for any given *y* in M, *y* satisfies [$H(a, y)$] classically.

Since $G(a, y)$ is recursive, it follows that there is some finite $k$ such that any Turing machine $T_1(y)$ that computes the recursive function $G(a, y)$ will halt and return the value 0 for $y = k$.

(*b*) Next, let $Q_2$ be the meta-assertion that $[H(a, y)]$ is effectively true in the standard interpretation M of PA, but that there is no uniformly effective method in M to determine that, for any given $y$ in M, $y$ satisfies $[H(a, y)]$ classically.

Since $G(a, y)$ is recursive, it follows that any Turing machine $T_2(y)$ that computes the arithmetical function $H(a, y)$ will halt, and return the symbol for self-termination (looping) for some value $y = k$.

(*c*) Finally, let $Q_3$ be the meta-assertion that $[H(a, y)]$ is effectively true in the standard interpretation M of PA, and that there is a uniformly effective method in M to determine that, for any given $y$ in M, $y$ satisfies $[H(a, y)]$ classically. We then have that that $[H(a, y)]$ is uniformly true in the standard interpretation M of PA.

Now, if we assume a Uniform Church Thesis, then $[H(a, y)]$ is PA-provable (since PA is, then, uniformly complete). Let $h$ be the Gödel-number of $[H(a, y)]$. We consider, then, Gödel's primitive recursive number-theoretic relation $xBy$[85], which holds in M if, and only if, $x$ is the Gödel-number of a proof sequence in PA for the PA-formula whose Gödel-number is $y$. It follows that there is some finite $k$ such that any Turing machine $T_3(y)$, which computes the characteristic function of $xBh$, will halt and return the value 0 for $x = k$.

---

[85] See footnote in Appendix, 1. Notation.

Since $Q_1$, $Q_2$, and $Q_3$ are mutually exclusive and exhaustive, it follows that, when run simultaneously over the sequence 1, 2, 3, ... of values for $y$, one of $T_1(y)$ / $T_2(y)$ / $T_3(y)$ will always halt for some finite value of $y$.

Thus, the halting problem is effectively solvable if we assume a Uniform Church Thesis. It, would then, also follow, first, that the parallel trio of Turing machines $\{T_1(y), T_2(y), T_3(y)\}$ is not a Turing machine; and, second, that the Turing thesis is false!

## II-5   Can all mental concepts be expressed mathematically?

### 5.1  What *is* mathematics?

Without attempting to address the issue in its broader dimensions, we may, without any loss of generality, consider mathematics simply as a set of precise, symbolic, languages.

Any language of such a set, say Peano Arithmetic PA (or Russell and Whitehead's Principia Mathematica, PM, or ZFC), is intended to express - in a finite, unambiguous, and communicable manner - relations between concepts that are external to the language PA (or to PM, or to ZFC). Each such language is, thus, in a Platonic sense, two-valued, if we assume that a relation either holds or does not hold externally (relative to the language).

Further, a selected, finite, number of primitive formal assertions about a finite set of selected primitive relations of, say, a language L are defined as axiomatically L-provable; all other assertions about relations that can be effectively defined in terms of the primitive relations are termed as L-provable if, and only if, there is a finite sequence of assertions of L, each of which is either a primitive assertion, or which can effectively be determined in a finite number of steps as an immediate consequence of any two assertions preceding it in the sequence by a finite set of rules of consequence.

An effective interpretation of a language L into another language, say PM (or PA, or ZFC, or English, etc.), is essentially the specification of an effective method by which any assertion of L is translated unambiguously into a unique assertion of PM (or PA, or ZFC, or English, etc.). Clearly, if a relation is provable in L, then it should be effectively decidable in any interpretation of L that shares a common logic - since a finite proof sequence of L would, prima facie, translate as a finite proof sequence in the interpretation.

The question arises: Is the converse true? In other words, if an assertion is decidable in an interpretation M of L, then does such decidability translate into an effective method of decidability in L?

Obviously, such a question can only be addressed unambiguously if there is an effective method for determining whether an assertion is decidable in M. If there is no such effective method, then we are faced with the following thesis:

**Thesis**: If there is no effective method for the unambiguous decidability of the assertions of a mathematical language L under an interpretation M, then L can only be considered a mathematical language of intuitive expression, but not a mathematical language of effective, and unambiguous, communication.

What this means is that, in the absence of an effective method of decidability in a mathematical language M that is a model of, say, a mathematical language such as PA, any correlation of *soundness*[86] between a PA-assertion and an assertion in M is essentially arguable; so it is meaningless to ask whether, in general, an assertion of PA is decidable under interpretation in M or not (the question of whether the assertion is decidable in PA or not is, then, an issue of secondary consequence).

---

[86] A formal language is, classically, said to be sound with respect to a given logic if every provable proposition of the language translates as a true sentence under any interpretation that shares the logic. The reverse condition is, classically, called completeness. (cf. also [WiSD]).

Our original question can, thus, be rephrased as:

> Can the Mind be considered as the standard - albeit intuitive - model M of PA, with the brain as its domain, and can PA be considered a formalisation of the Mind?

The first part of the question can, of course, be easily answered affirmatively, since every proof sequence of PA can, intuitively, be recognised by, and thus be seen as an effective method of decidability in, M.

The second part, on the other hand, can be answered negatively just as easily, if we accept that the Mind can be seen to experience hallucinations (reflecting synaptic patterns that, by definition, we assume correspond to some aspects of the physical state of the brain) that cannot conceivably be effectively verified in any formal language.

However, if we adopt the Individual and Uniform Church Theses §I-5(*vii*)-(*viii*) above, then every individually effective method, and every uniformly effective method of decidability, or computability, in the Mind corresponds to some effective method of decidability, or computability, in the language that is used to describe the Mind.

It follows that, if any two interpretations of such a language are isomorphic (as indicated in §I-6(*xiv*) above), then the language can, in a sense, be taken as adequately formalising that part of the activity of the brain that lends itself to mathematical representation.

## 5.2 Knowledge and intuitive truth

The alternative consequences of Godel's reasoning discussed in these essays draw significantly upon the way we choose to perceive the nature of *Intuitive Knowledge*, and more particularly the nature of factual, or intuitive truth.

By *Intuitive Knowledge* we refer loosely to that body of proactive knowledge that stems directly from our conscious states, in contrast to our reactive *Instinctive Knowledge*, which stems from, and lies within, our sub-conscious and unconscious states.

These essays are intended to highlight the wider significance of an issue that may otherwise lie in obscurity due to the specialised nature of the subject - that the influence, on our current modes of thought, of the interpretations, and conclusions, drawn from Gödel's original paper [Go31a] may have a wider, multi-disciplinary, element that is not obvious from an appreciation of its purely logical and mathematical import.

## 5.3 Implicit influence of Gödel's Platonism

The roots of this influence may be traced to implicitly Platonistic elements that underlie classical first order Peano Arithmetic, PA, which is based essentially on the formal system P defined by Gödel in his paper [Go31a]. Loosely speaking Gödel, who was an explicitly strong Platonist, assumed the existence of a world of ideals that could objectively be referred to for arbitrating which of our assumptions or premises, when formally expressed in a rigorously constructed scientific language, were *intuitively true*, and which were *intuitively false*, under any given interpretation.

Now one may, when attempting to express mental concepts within a language, argue reasonably - as Gödel does - in Platonistic terms, and define *intuitive truths* as characteristics of relationships that are assumed to exist in some absolute sense (that is, even in the absence of any perceiver) between the objects of an external ontology (both of which are also taken to exist in some absolute sense).

We argue, however, that, for effective communication, we may need an alternative view of relationships as belonging to individual perceptions that we consciously construct, and selectively assign, to abstract objects (that themselves are individual conceptual constructs) of an abstract ontology (that is similarly an individual conceptual construct).

In other words, each individual perception can, reasonably, be assumed to be a subjective, abstract, construct, based on a unique, one-of-a-kind, never-to-be-repeated, consciousness of an individual experience. An intuitive truth is, then, essentially, a constructed, space-time localised, individual *factual truth* (we shall, henceforth, use the terms synonymously). It corresponds to a subjectively constructed characteristic of the expression of an individual perception. Loosely speaking, it corresponds to a characteristic of the way we construct an expression for that which we select as common to a series of subjective perceptions, rather than to a characteristic that we discover of an objectively observed *something*.

The distinction seems significant. Platonic concepts, when introduced into the interpretations of a formal language of communication, could permit us to misleadingly validate our subjective intuitive perceptions of individual factual truth as being reflective of some absolute Truth that must be of universal significance in a Utopian, Platonistic world. The Biblical Tower of Babel can, arguably, be seen as illustrating the extreme possibilities of such Platonistic beliefs.

## 5.4 Formal algorithmic and non-algorithmic truth

Now, a constructive interpretation of Gödel's Theorems is that they actually establish the effectiveness of our ability to also communicate abstractions that we intellectually conceive as non-algorithmic[87], on the basis of our individual sensory perceptions. They thus have to do with the efficiency and effectiveness of our language of communication, and involve the concept of logical, or formal, truth (we treat the two terms as synonymous).

---

[87] We contrast this with a classical interpretation of Gödel's Theorems such as is reflected in Davis' remark [Da90], "There is certainly room for disagreement about whether the processes by which mathematical (or physical) theories are developed and accepted are algorithmic. But Gödel's theorem has nothing decisive to contribute to the discussion."

We argue that, in formal languages, a selected set of axiomatic truths i s expressed as a set of *Axioms* (*or Axiom schemas*). The selection criteria is that the Axioms are readily accepted by any perceiver as faithfully reflecting some significant factual truths, pertaining to the expression of abstract constructed elements of the constructed ontology under consideration, as perceived and conceived by the perceiver.

For the most basic, and intuitive, of our scientific languages, namely Number Theory or our Arithmetic of the natural numbers, we take the commonly accepted selection of such axiomatic truths as the classical set of Peano's Postulates, first expressed in semi-axiomatic format by Dedekind in 1901.

The challenge, then, is to express these in a formal language such as Peano Arithmetic, PA, along with a suitable set of *Rules of Inference* by which we can assign unique formal truth-values algorithmically to as many well-formed propositions of the theory as possible that are not Axioms.

The concept of *formal algorithmic truths* is, thus, merely the result of the application of a set of Rules of Inference for effectively assigning such formal truth-values algorithmically to various logical permutations and combinations of axiomatic truths as (finite and infinite) compound assertions (which, ideally, should not introduce any new axiomatic truths that are not already implicit within the Axioms).

Now, Tarski's Theorem, based on Gödel's remarkable reasoning [Go31a], is that, in any such system of Arithmetic, there are also formal truths that cannot be verified algorithmically. In other words, there are *formal non-algorithmic truths* that are consequences of the axiom system, but which can only be verified non-algorithmically.

By choosing to adopt Tarski's definitions of satisfiability and truth without effective verifiability, classical interpretations of formal, non-algorithmic, truth, hold that such truth is a property of the interpretation of a language, and so it cannot be contained within the language.

The limiting consequences of such interpretations on the beliefs of other disciplines, which look to mathematics for providing them with languages for the pr ecise expression, and unambiguous communication, of scientific concepts, is seen in Gurney's following remarks [Gu96], in the context of neural nets:

"... contemporary 'AI engines' are still vehicles for the instantiation of the theoretic stance which claims that cognition can be described completely as a process of formal, algorithmic symbol manipulation. ...

AI has not fulfilled much of the early promise that was conjectured by the pioneers in the field. ... Principal among these are the belief that all knowledge or information can be formalised, and that the mind can be viewed as a device that operates on information according to formal rules. It is precisely in those domains of experience where it has proved extremely difficult to formalise the environment, that the 'brittle' rule-based procedures of AI have failed.

The differentiation of knowledge into that which can be treated formally and that which cannot ... makes the distinction between cultural, or public knowledge and private, or intuitive knowledge. The stereotypical examples of the former are found in science and mathematics, whereas the latter describes, for instance, the skills of a native speaker or the intuitive knowledge of an expert in some field. In the connectionist view, intuitive knowledge cannot be captured in a set of linguistically formalized rules and a completely different strategy must be adopted."

In these essays, we argue, however, that such a constricted vision, of the range, and capability, of intellectual pursuit, as is implicit in the above remarks, may be unnecessary. A constructive interpretation of formal non-algorithmic truth may, indeed, be able to introduce effective methods of verifiability into Tarski's definitions, and to extend Church's Thesis, so that such verifiability follows from the language itself, and is independent of any interpretation.

Such a language could, then, be defined as effectively, if not uniformly, complete.

## 5.5 Significance of formal truth

Now, the significance of formal truths lies in our experience that the individual factual truths of our perceptions can generally be corresponded in a communicable language with a high degree of correlation to the formal truths of the language.

This appears to suggest that the evolutionary significance to us (and possibly to any intelligence whose evolution is based on communication) of any set of individual factual truths may be proportional to the body of formal truths that can be constructed by various logical combinations and permutations of the axiomatic truths that we collectively accept as representative of some of our individual factual truths.

## 5.6 Non-verifiable interpretations of Gödel's formal system

The significance of Gödel's Theorems lies in the fact that they are derived in a system of Axioms where the Rules of Inference lead to a particularly rich body of expressions that can be assigned formal truth-values under various interpretations of the symbols of the theory. However, a major feature of such a system is that it also lends itself to interpretations of the chosen Rules of Inference that are non-constructive, in the sense that they are able to assign, implicitly and sweepingly, non-verifiable formal truth-values in some models to various expressions. Thus the language, in a sense, admits formally true expressions under some interpretations that cannot be correlated, even in principle, to any factual truths of a human perception.

However, this is not necessarily the drawback that it appears to be at first sight. In fact such duality should both encourage and discipline us in our use of rich and creative languages that can be used, both, as languages for the mathematical expression of some of our mental concepts, and as languages of effective communication for some of our mental concepts. This should, moreover, force us to focus on devising criteria by which we can recognise vague, or ill-defined, concepts that, although expressible mathematically, cannot be effectively communicated by, and hence cannot be introduced as primitive symbols of (or defined in terms of existing primitive symbols of), a formal language.

## II-6 Can a constructive interpretation of Peano Arithmetic model some of the more paradoxical concepts of Quantum Mechanics?

Interestingly, insistence on a constructive verifiability of Tarski's definitions of the satisfiability and truth of the formulas of a formal language under an interpretation, as suggested in this essay, implies that PA can also express relations that are deterministic, yet essentially unpredictable; such a language would have significance for the expression of natural phenomena that are best described in quantum mechanical terms.

Thus, in Anand [An03i], we indicate how the introduction of constructive definitions of classical mathematical concepts permits formal systems of Peano Arithmetic to model some of the more paradoxical concepts of Quantum Mechanics. For instance, as a consequence of §I-6(*xiii*), consider the following argument:

(*a*) Gödel has proved, in his 1931 paper [G031a], that there is an arithmetic formula [*R*(*x*)] such that, for any given *k*, [*R*(*k*)] is provable.

(*b*) Hence, for any given *k*, there is always some effective method for evaluating the arithmetic expression *R*(*k*) (when treated as a Boolean function).

(*c*) Gödel has also proved in the above paper that [(A*x*)*R*(*x*)] is not provable.

(*d*) *Thesis*: There is no uniformly effective method (algorithm/Turing machine) that can evaluate the arithmetic expression *R*(*x*) for any given *x*.

(*e*) Thus, *R*(*n*) is individually computable, but not uniformly computable.

(*f*) *Theorem (provable by induction)*: For any given *k*, we can always find some effective method (algorithm/Turing machine) *T*(*k*) that can compute *R*(*n*) for all *n*<*k*, i.e. *T*(*k*) terminates for all *n*<*k*, but it "loops" on input *k*. (Note: All methods that evaluate *R*(*n*) for all *n*<*k* cannot be non-terminating on input *k*; this would imply that *R*(*k*) is undefined, which would contradict (*b*).)

(*g*) *Quantum interpretation*: The process of finding *T*(*k*+1) can be corresponded, first, to the act of finding a suitable method of measuring the value *R*(*k*) precisely, and, second, to the collapse of the wave function at *k* as a result of the measurement; we then have the new "state" *T*(*k'*), which can evaluate the value of *R*(*n*) for all *n*<*k'*, where *k*<*k'*, but not beyond!

(*h*) If, now, we have some law that determines the state *T*(*k'*) from the state *T*(*k*) and the interaction at *k*, we have a deterministic interaction that is, nevertheless, absolutely unpredictable, where we may then define free will as absolute unpredictability. (We note that, if *k'* > *k*+1, we have a language that admits inter-actions that can leave the state *T*(*k'*) unchanged.)

Now we note that a counter-thesis to (*d*) would be:

(*i*) *Counter-Thesis*: There is some uniformly effective method (algorithm/Turing machine) that can evaluate the arithmetic expression $R(x)$ for any given $x$.

It follows from Gödel's reasoning in [Go31a] that both (*d*) and (*i*) are effectively unverifiable, since they cannot be proved formally. We thus have two standard models of Peano Arithmetic - classical and constructive - that are mutually inconsistent. If we assume that both are consistent, the above argument indicates the interpretation that implies (*d*) may be the more suitable language for expressing, and effectively communicating, some of the more paradoxical concepts of classical Quantum Theory.

## III    In conclusion

In this essay, we have argued that classical interpretations of Gödel's, Tarski's, Turing's and Cantor's formal reasoning, and of their conclusions, implicitly imply that mathematical languages are essentially incomplete, in the sense that the truth of some arithmetical propositions of any formal mathematical language, under any interpretation, is, both, non-algorithmic, and essentially unverifiable.

We have argued, however, that a language of general, scientific, discourse, which intends to mathematically express, and unambiguously communicate, intuitive concepts that correspond to scientific views of the universe, cannot allow its mathematical propositions to be interpreted ambiguously. Such a language must, therefore, define mathematical truth verifiably.

We have considered a constructive interpretation of classical, Tarskian, truth, and of Gödel's reasoning, under which any formal system of Peano Arithmetic - classically accepted as the foundation of all our mathematical languages - is verifiably complete in the above sense.

We have, further, indicated how some of the more contentious foundational issues regarding the nature of mathematical objects and mathematical truth, and some seemingly paradoxical mathematical concepts - such as that of Gödel's undecidable proposition, that of Turing's uncomputable Halting function, and that of a probabilistic, yet deterministic, Quantum mechanics - can be expressed, and interpreted, naturally under a constructive definition of mathematical truth.

We have, thus, argued that there are no absolute limits on our capacity to express human cognition unambiguously in a communicable language; there are only temporal limits - not necessarily absolute - to the capacity of classical interpretations to communicate unambiguously that which we intended to capture within our formal expression.

## Appendix

### 1. Notation

Unless specified otherwise, we generally follow the notation introduced by Mendelson in his English translation of Gödel's 1931 paper [Go31a]; however, for convenience of exposition, we refer to it as Gödel's notation. Two notable exceptions: we use the notation "$(Ax)$", whose classical, standard, interpretation is "for all $x$", to denote Gödel's special symbol for Generalisation; the successor symbol is denoted by "$S$", instead of by "$f$".

Following Gödel (cf. [Go31a], footnote 13), we use square brackets to indicate that the expression [(A$x$)], including square brackets, only denotes the uninterpreted string[88] named[89] within the brackets. Thus, [(A$x$)] is not part of the formal system P, and would be replaced by Gödel's special symbolism for Generalisation wherever it occurs.

Following Gödel's definitions of well-formed formulas[90], we note that juxtaposing the string [(A$x$)] and the formula[91] [$F(x)$] is the formula [(A$x$)$F(x)$], juxtaposing the symbol [~] and the formula [$F$] is the formula [~$F$], and juxtaposing the symbol [v] between the formulas [$F$] and [$G$] is the formula [$F$v$G$].

The number-theoretic functions and relations are defined explicitly by Gödel[92] [Go31a]. The formulas are defined implicitly by his reasoning.

## 2. Some of Gödel's definitions

For convenience, we reproduce some of the more significant of Gödel's definitions. We take P to be Gödel's formal system[93], and define ([Go31a], Theorem VI, p24-25):

---

[88] We define a "string" as any concatenation of a finite set of the primitive symbols of the formal system under consideration.

[89] We note that the "name" inside the square brackets only serves as an abbreviation for some string in P.

[90] We note that all well-formed formulas of P are strings of P, but all strings of P are not well-formed formulas of P.

[91] By "formula", we shall henceforth mean a "well-formed formula" as defined by Gödel ([Go31a], p11).

[92] We follow Gödel's definition of recursive number-theoretic functions and relations ([Go31a], p14-17). We note, in particular, that Gödel's recursive number-theoretic function $Sb(x \; 19|Z(y))$ is defined as the Gödel-number of the P-formula that is obtained from the P-formula whose Gödel-number is $x$ by substituting the numeral [$y$], whose Gödel-number is $Z(y)$, for the variable whose Gödel-number is 19 wherever the latter occurs free in the P-formula whose Gödel-number is $x$ ([Go31a], p20, def.31). We also note that Gödel's recursive number-theoretic relation $xBy$ holds if, and only if, $x$ is the Gödel-number of a proof sequence for the P-formula whose Gödel-number is $y$ ([Go31a], p22, def. 45).

[93] Gödel ([Go31a], p9-13).

(*i*)     $Q(x, y)$ as Gödel's recursive number-theoretical relation $\sim xB(Sb(y\ 19|Z(y)))$.

(*ii*)    $[R(x, y)]$ as a formula that represents $Q(x, y)$ in the formal system P.

(*iii*)   $q$ as the Gödel-number[94] of the formula $[R(x, y)]$ of P.

(*iv*)    $p$ as the Gödel-number of the formula $[(Ax)R(x, y)]$[95] of P.

(*v*)     $[p]$ as the numeral that represents the natural number $p$ in P.

(*vi*)    $r$ as the Gödel-number of the formula $[R(x, p)]$ of P.

(*vii*)   $17Genr$ as the Gödel-number of the formula $[(Ax)R(x, p)]$ of P.

(*viii*)  $Neg(17Genr)$[96] as the Gödel-number of the formula $[\sim(Ax)R(x, p)]$ of P.

(*ix*)    $R(x, y)$ as the standard interpretation of the formula $[R(x, y)]$ of P.

(*x*)     $Wid$(P) as the number-theoretic assertion $(Ex)(Form(x)\ \&\ \sim Bew(x))$[97].

          (We note that $Wid$(P) is defined by Gödel ([Go31a], p36) as equivalent to the
          meta-assertion "P is consistent".)

(*xi*)    $[Con(P)]$ as the formula that represents $Wid$(P) in the formal system P.

---

[94] By the "Gödel-number" of a formula of P, we mean the natural number corresponding to the formula in the 1-1 correspondence defined by Gödel ([Go31a], p13).

[95] We note that "[(Ax)][R(x, y)]" and "[(Ax)R(x, y)]" denote the same formula of P.

[96] We note that Gödel's recursive number-theoretic function $Neg(x)$ is the Gödel-number of the P-formula that is the negation of the P-formula whose Gödel-number is $x$ ([Go31a], p18, def. 13).

[97] We note that Gödel's recursive number-theoretic relation $Form(x)$ is satisfied if, and only if, $x$ is the Gödel-number of a P-formula ([Go31a], p19, def. 23). Also, Gödel's number-theoretic relation $Bew(x)$ is satisfied if, and only if, $x$ is the Gödel-number of a provable P-formula ([Go31a], p22, def. 46).

(*xii*)  *w* as the Gödel-number of the formula [*Con*(P)] of P [1, p37].

(*xiii*) *Con*(P) as the standard interpretation of the formula [*Con*(P)] of P.

## References[98]

[An02a]  Anand, B. S. 2002. *Reviewing Gödel's and Rosser's meta-reasoning of undecidability.* (*Web essay*)

<*Preprint*: http://alixcomsi.com/Constructivity_consider.htm>

[An02b]  Anand, B. S. 2002. *Omega-inconsistency in Gödel's formal system: a constructive proof of the Entscheidungsproblem.* (*Web essay*)

<*Preprint*:  http://alixcomsi.com/CTG_00.htm>

[An02c]  Anand, B. S. 2002.  *Some consequences of defining mathematical objects constructively and mathematical truth effectively.* (*Web essay*)

<*Preprint*: http://alixcomsi.com/CTG_06_Consequences.htm>

[An03a]  Anand, B. S. 2003. *Is there a duality in the classical acceptance of non-constructive, foundational, concepts as axiomatic?* (*Web essay*)

<*Preprint*: http://alixcomsi.com/CTG_06_Consequences_Bringsjord.htm>

[An03b]  Anand, B. S. 2003. *Three beliefs that lend illusory legitimacy to Cantor's diagonal argument?* (*Web essay*)

<*Preprint*: http://alixcomsi.com/Three_beliefs.htm>

---

[98] This essay has liberally used the cited on-line resources, not necessarily peer-reviewed, and others that are freely available on the world-wide web.

[An03c]  Anand, B. S. 2003. *Is there a "loophole" in Gödel's interpretation of his formal reasoning and its consequences?* (*Web essay*)

    *<Preprint*: http://alixcomsi.com/Is_there_a_loophole.htm>

[An03d]  Anand, B. S. 2003. *Can Turing machines capture everything we can compute?* (*Web essay*)

    *<Preprint*: http://alixcomsi.com/Can_Turing_machines.htm >

[An03e]  Anand, B. S. 2003. *The formal roots of Platonism* (*Web essay*)

    *<Preprint*: http://alixcomsi.com/The_formal_roots_of_Platonism.htm>

[An03f]  Anand, B. S. 2003. *Can we express every transfinite concept constructively?* (*Web essay*)

    *<Preprint*: http://alixcomsi.com/Can_we_express_every_transfinite.htm>

[An03g]  Anand, B. S. 2003. *Is the Halting probability a Dedekind real number?* (*Web essay*)

    *<Preprint*: http://alixcomsi.com/Is_the_Halting_probability.htm>

[An03h]  Anand, B. S. 2003. *Why we must heed Wittgenstein's "notorious paragraph".* (*Web essay*)

    *<Preprint*: http://alixcomsi.com/Why_we_must_heed.htm>

[An03i]  Anand, B. S. 2003. *How definitive is the standard interpretation of Gödel's Incompleteness Theorem?* (*Web essay*)

    *<Preprint*: http://alixcomsi.com/How_definitive_is_the_standard.htm>

[An03j] Anand, B. S. 2003. *Can Laplace's formula model a deterministic universe that is irreducibly probabilistic?* (*Web essay*)

*Preprint*: http://alixcomsi.com/Can_Laplace's_formula.htm>

[An04a] Anand, B. S. 2004. *How definitive is the standard interpretation of Goodstein's argument?* (*Web essay*)

*Preprint*: http://alixcomsi.com/Goodstein_argument.htm>

[Br93] Bringsjord, S. 1993. *The Narrational Case Against Church's Thesis.* Easter APA meetings, Atlanta. (*Web essay*)

*Web-page*: http://www.rpi.edu/~brings/SELPAP/CT/ct/ct.html>

[Da90] Davis, M. 1990. *Is Mathematical Insight Algorithmic?* Behavioural and Brain Sciences, vol. 13 (1990), pp. 659-660.

*Davis' web-page*: http://www.cs.nyu.edu/cs/faculty/davism/>

[DS99] Dehaene, S. Spelke, E. Pinet, P. Stanescu, R. Tsivkin S. 1999. *Sources of mathematical thinking: behavioral and brain-imaging evidence.* Science, 7 May 1999, vol. 284, 970{974.

[Go31a] Gödel, Kurt. 1931. *On formally undecidable propositions of Principia Mathematica and related systems I.* Translated by Elliott Mendelson. In M. Davis (ed.). 1965. The Undecidable. Raven Press, New York.

[Go31b] Gödel, Kurt. 1931. *On formally undecidable propositions of Principia Mathematica and related systems I.* Translated by B. Meltzer. 1962. University of Edinburgh.

*Web page*: http://home.ddc.net/ygg/etext/godel/index.htm>

[Go51]   Gödel, Kurt. 1951. *Some basic theorems on the foundations of mathematics and their implications.* Gibbs lecture. In Kurt Gödel, Collected Works III, pp. 304-323. 1995. Unpublished Essays and Lectures. S. Feferman et al (ed.). Oxford University Press, New York.

[Gu96]   Gurney, K. 1996. *Neural Nets.* Notes to "An Introduction to Neural Networks". UCL Press, London. (*Web essay*)

   *<Web page*: http://www.shef.ac.uk/psychology/gurney/notes/l10/subsection3_1_1.html>

[Gw02]   Gowers, W. T. 2002. *Does mathematics need a philosophy?* Presented before the Cambridge University Society for the Philosophy of Mathematics and Mathematical Sciences, Michaelmas 2002. (*Web essay*)

   *<Web page*: http://www.dpmms.cam.ac.uk/~wtg10/philosophy.html>

[Ha47]   Hardy, G.H. 1947, 9th ed. Pure Mathematics. Cambridge, New York.

[He04]   Christer-Hennix, C. 2004. *Some remarks on Finitistic Model Theory, Ultra-Intuitionism and the main problem of the Foundation of Mathematics.* ILLC Seminar, 2nd April 2004, Amsterdam.

[HA00]   Hodges, A. 2000. *Uncomputability in the work of Alan Turing and Roger Penrose.* Talk given for Interface 5, Hamburg, 6 October 2000

   *<Web page*: http://www.turing.org.uk/philosophy/lecture1.html>

[HW98]  Hodges, W. 1998. *An editor recalls some hopeless papers.* The Bulletin of Symbolic Logic, Volume 4, Issue 1, March 1998, pages 1-16.

   *<Web link*: http://plato.stanford.edu/archives/win2001/entries/tarski-truth/>

[HW01]  Hodges, W. 2001. *Tarski's Truth Definitions.* Stanford Encyclopedia of Philosophy. (Winter 2001 Edition). Edward N. Zalta (ed.). (*Web essay*)

       *<Web page*: http://plato.stanford.edu/archives/win2001/entries/tarski-truth/>

[Ka59]  Kalmár, L. 1959. *An Argument Against the Plausibility of Church's Thesis.* In Heyting, A. (ed.) Constructivity in Mathematics. North-Holland, Amsterdam.

[Kl36]  Kleene, S.C. 1936. *General Recursive Functions of Natural Numbers.* Math. Annalen **112**.

[La51]  Landau, E.G.H. 1951. Foundations of Analysis. Chelsea Publishing Co., New York.

[Lu61]  Lucas, J. R. 1961. *Minds, Machines and Gödel.* Philosophy, XXXVI, 1961, pp.(112)-(127); reprinted in The Modeling of Mind, Kenneth M.Sayre and Frederick J.Crosson, eds., Notre Dame Press, 1963, pp.[269]-[270]; and Minds and Machines, ed. Alan Ross Anderson, Prentice-Hall, 1954, pp.{43}-{59}

       *<Web page*: http://users.ox.ac.uk/~jrlucas/Godel/mmg.html>

[Ma02]  Manin, Yu. I. 2002. *Georg Cantor and his heritage.* Meeting of the German Mathematical Society and the Cantor Medal award ceremony. (*Web essay*)

       *<Web page*: http://arxiv.org/PS_cache/math/pdf/0209/0209244.pdf>

[Me64]  Mendelson, Elliott. 1964. Introduction to Mathematical Logic. Van Norstrand, Princeton.

[Me90]  Mendelson, E. 1990. *Second Thoughts About Church's Thesis and Mathematical Proofs.* Journal of Philosophy **87.5**.

[Pe90]   Penrose, R. (1990, Vintage edition). The Emperor's New Mind: Concerning
         Computers, Minds and the Laws of Physics. Oxford University Press.

[Pe94]   Penrose, R. (1994). Shadows of the Mind: A Search for the Missing Science of
         Consciousness. Oxford University Press.

[Po01]   Podnieks, Karlis. 2001. What is Mathematics: Goedel's Theorem and Around.

         <*e-textbook*: http://www.ltn.lv/~podnieks/gt.html>

[RH01]   Ramachandran, V. S. and Hubbard, E. M. 2001. *Synaesthesia - A Window Into
         Perception, Thought and Language.* Journal of Consciousness Studies, 8, No.
         12.

         <*Ramachandran's Web page*: http://psy.ucsd.edu/chip/ramapubs.html>

[Ru53]   Rudin, Walter. 1953. Principles of Mathematical Analysis. McGraw Hill, New
         York.

[Rs01]   Russell, B. A. W. 1901. *Recent Work on the Principles of Mathematics.*
         International Monthly, 4, 83-101. Repr. as *Mathematics and the Metaphysicians*
         in Russell, Bertrand, *Mysticism and Logic*, London: Longmans Green, 1918, 74-
         96.

[Sh91]   Shelah, S. 1991. *The Future of Set Theory*. Proceedings of the Bar Ilan Winter
         School (January 1991). Israel Mathematical Conference proceedings, vol 6.

         <*Web page*: http://shelah.logic.at/E16/E16.html>

[Ta36]   Tarski, A. 1936. *Der Wahrheitsbegriff in den formalisierten Sprache.* Studia
         Philos., Vol. 1. Expanded English translation, *The concept of truth in the
         languages of the deductive sciences*, in *Logic, Semantics, Metamathematics,*

*papers from 1923 to 1938* (p152-278), ed. John Corcoran. 1983. Hackett
Publishing Company, Indianapolis.

[Ti61]   Titchmarsh, E. C. 1961. The Theory of Functions. Oxford University Press.

[Tu36]   Turing, Alan. 1936. *On computable numbers, with an application to the
Entscheidungsproblem.* Proceedings of the London Mathematical Society, ser.
2. vol. 42 (1936-7), pp.230-265; corrections, Ibid, vol 43 (1937) pp. 544546.

<*Web page*: http://www.abelard.org/turpap2/tp2-ie.asp - index>

[Wi78]   Wittgenstein, Ludwig. (1978 edition). Remarks on the Foundations of
Mathematics. MIT Press, Cambridge.

[WiPM]   *Philosophy of Mathematics.* Wikipedia. (*Web essay*)

<*Web page*: http://en.wikipedia.org/wiki/Philosophy_of_mathematics>

[WiSD]   *Soundness.* Wikipedia. (*Web essay*)

<*Web page*: http://en.wikipedia.org/wiki/Soundness>

(*Updated: Monday 5th July 2004 7:50:11 AM by re@alixcomsi.com*)